\documentclass[11pt]{article}
\usepackage{amsmath,amssymb}
\usepackage[latin1]{inputenc}
\usepackage[dvips]{graphicx}

\topskip=\baselineskip  \headsep=0in
\newtheorem{trm}{Theorem}[section]
\newtheorem{lem}[trm]{Lemma}
\newtheorem{prop}[trm]{Proposition}

\newtheorem{remark}[trm]{Remark}
\newcommand{\qed}{\enspace\vrule height6pt width4pt depth2pt}

\newenvironment{proof}{\par\noindent{\bf Proof.}}{$\qed$\par\bigskip}

\newcommand{\GEN}[1]{\langle #1 \rangle}

\newcommand{\U}{{\cal U}}

\newcommand{\vp}{\varphi}
\newcommand{\Char}{{\rm char}}
\newcommand{\z}{\ensuremath{\mathcal{Z}}}

\date{}
\title{Antisymmetric elements in group rings with an orientation morphism\footnote{The first author has been
partially supported by FAEPEMIG of Brazil. The
second author has been partially supported by
Onderzoeksraad of Vrije Universiteit Brussel, Fonds
voor Wetenschappelijk Onderzoek (Belgium),
Flemish-Polish bilateral agreement BIL2005/VUB/06.
The third author has been partially supported by
D.G.I. of Spain and Fundaci\'on S\'eneca of Murcia.
\newline 2000 Mathematics Subject Classification: 16S34, 16W10, 16U60.}}
\author{O. Broche \and E. Jespers \and M. Ruiz}

\begin{document}

\setcounter{page}{0} \thispagestyle{empty}

\vspace{10cm}
\begin{center}
{\LARGE Antisymmetric elements in group rings with
an orientation morphism}

\vspace{5cm}

\begin{tabular}{lll}
O. Broche Cristo &Eric Jespers &  Manuel Ruiz
\\ Dep. de Ci\^{e}ncias Exatas &Dept. Mathematics  & Dep.
M\'{e}todos Cuantitativos e Inform\'{a}ticos \\
Universidade Federal de Lavras &Vrije Universiteit
Brussel &Universidad Polit\'{e}cnica de Cartagena
\\ Caixa Postal 3037& Pleinlaan 2 &  Paseo Alfonso
XIII, 50 \\ 37200-000 Lavras, Brazil  & 1050
Brussel, Belgium
 & 30.203 Cartagena, Spain \\ osnel@ufla.br & efjesper@vub.ac.be
  & Manuel.Ruiz@upct.es
\end{tabular}

\vspace{2cm}

{\bf Proposed Running title}: Antisymmetric elements
in group rings

\end{center}
\newpage

\maketitle

\begin{abstract}
Let $R$ be a commutative ring, $G$ a group and $RG$
its group ring. Let $\varphi_{\sigma} :
RG\rightarrow RG$ denote the involution defined by
$\varphi_{\sigma} (\sum r_{g}g) = \sum r_{g} \sigma
(g) g^{-1}$, where $\sigma:G\rightarrow \{\pm 1\}$
is a group homomorphism  (called an orientation
morphism). An element $x$ in $RG$ is said to be
antisymmetric if $\varphi_{\sigma} (x) =-x$. We give
a full characterization of the groups $G$ and its
orientations for which the  antisymmetric elements
of $RG$ commute.
\end{abstract}

\section{Introduction}
Let $R$ be a commutative ring with unity and let $G$
be a group. Let $\varphi$ be an involution on the
group ring $RG$. Denote by $\U(RG)$ the group of
units of the group ring $RG$ and  by
$(RG)^{-}_{\vp}$ the set of  its antisymmetric
elements, that is, $$(RG)^{-}_{\vp}=\{\alpha\in RG
\mid \vp(\alpha)=-\alpha\}.$$
In this paper we investigate when $(RG)^{-}_{\vp}$
is commutative, that is $ab=ba$ for all $a,b\in
(RG)^{-}_{\vp}$.

The group of $\varphi$-unitary units of $RG$ is
defined by $$\mathcal{U}_{\varphi} (RG)=\{
u\in\U(RG) \mid u\varphi (u)=1\}.$$ For general
algebras there is a close relationship between the
$\varphi$-unitary units and the antisymmetric
elements. For example, in \cite{GM} Giambruno and
Polcino Milies show that if $\varphi$ is an
involution on a finite dimensional semisimple
algebra $A$ over an algebraically closed field $F$
with $char(F)\neq 2$ then $\mathcal{U}_{\varphi}(A)$
satisfies a group identity if and only if
$(A)^{-}_{\varphi}$ is commutative. Moreover, if $F$
is a nonabsolute field then
$\mathcal{U}_{\varphi}(A)$ does not contain a free
group of rank $2$ if and only if $(A)^{-}_{\varphi}$
is commutative. Giambruno and Sehgal, in
\cite{gia-seh}, showed that if $B$ is a semiprime
ring with involution $\varphi$, $B=2B$ and
$(B)^{-}_{\varphi}$ is Lie nilpotent then
$(B)^{-}_{\varphi}$ is commutative and $B$ satisfies
a polynomial identity of degree $4$. The latter
shows that crucial information of the algebraic
structure of $A$ can be determined by that of
$(A)^-_\vp$. We state two more important results of
this nature.
 Amitsur in
\cite{A} proves that for an arbitrary algebra $A$
with an involution $\varphi$, if $A_{\varphi}^{-}$
satisfies a polynomial identity (in particular when
$A_{\varphi}^{-}$ is commutative) then $A$ satisfies
a polynomial identity. Gupta and Levin in
\cite{gupta} proved that for all $n\geq 1$
$\gamma_n(\U(A))\leq 1+ L_n(A)$. Here $\gamma_n(G)$
denotes the $n$th term in the lower central series
of the group $G$ and $L_n(A)$ denotes the two sided
ideal of $A$ generated by all Lie elements of the
form $[a_1,a_2,\dots,a_n]$ with $a_i\in A$ and
$[a_1]=a_1$, $[a_1,a_2]=a_1a_2-a_2a_1$ and
inductively
$[a_1,a_2,\dots,a_n]=[[a_1,a_2,\dots,a_{n-1}],a_n]$.
Smirnov and Zalesskii in \cite{ZS}, proved that, for
example,
 if the Lie ring generated by the elements of the form $g+g^{-1}$
 with $g\in \U(A)$ is Lie nilpotent then $A$ is Lie nilpotent.

 Special
attention has been given to the classical involution
$*$ on $RG$ that is the $R$-linear map defined by
mapping $g\in G$ onto $g^{-1}$. In case $R$ is a
field of characteristic $0$ and $G$ is a periodic
group, Giambruno and Polcino Milies in \cite{GM}
described when $\mathcal{U}_{*}(RG)$ satisfies a
group identity. Gon\c{c}alves and Passman in
\cite{GoPas} characterized when
$\mathcal{U}_{*}(RG)$ does not contain non abelian
free groups when $G$ is a finite group and $R$ is a
nonabsolute field. Giambruno and Sehgal, in
\cite{gia-seh}, show that if $R$ is a field of
characteristic $p\geq 0$, with $p\neq 2$ and $G$  a
group without $2$-elements, then the Lie nilpotence
of $(RG)_{*}^{-}$ implies the Lie nilpotence of
$RG$. Giambruno, Polcino Milies  and Sehgal in
\cite{GM,gs} characterized when $(RG)_{*}^{-}$ is
Lie nilpotent.

Because of all the above mentioned results, it is
relevant to determine when the antisymmetric
elements of a group ring commute. Recently, for an
arbitrary involution $\vp$ on a group
 $G$ (extended by linearity to $RG$) and a commutative ring $R$, Jespers and Ruiz \cite{JR2}
characterized  when  $(RG)^{-}_{\vp}$ is
commutative.
 This generalizes earlier work of Broche and Polcino Milies \cite{BP} in case $\vp$ is the classical involution.
 The characterizations obtained in both papers are in terms of the algebraic structure of some subgroups of $G$.

In \cite{osnel2}, \cite{bovdi}, \cite{bovdisehg}
and \cite{li2} various authors considered
involutions on a group ring $RG$ that are not
determined by $R$-linearity by an
 involution on $G$.
The following is an example of such an involution
$\varphi_{\sigma}$ that was introduced  by Novikov
in \cite{novikov} in the context of K-theory and
algebraic topology:
\begin{equation*}\label{novi}
\varphi_{\sigma} \left(\sum\limits_{g\in G}\alpha_g
g\right)= \sum\limits_{g\in G}
\alpha_g\sigma(g)g^{-1},
\end{equation*}
where $\sigma:G\rightarrow \{\pm 1\}$ is a group
homomorphism (called  an orientation of $G$) and all
$\alpha_{g}\in R$. Note that such a $\sigma$ is
uniquely determined by its kernel $\ker (\sigma )
=N$.

The aim of this paper is to prove the following
theorem in which we fully describe when
$(RG)_{\varphi_{\sigma}}^{-}$ is commutative, and
this in terms of presentations of the groups $G$ and
kernels $N$. Because of the results mentioned above,
we will only deal with the case that   $G\neq \ker
(\sigma )$ and therefore $\Char(R)\neq 2$. Moreover,
if  $\Char(R)=2$ then the antisymmetric elements are
precisely the symmetric elements and in
\cite{osnel2,JR} it has been classified when the
symmetric elements in $RG$ commute.

We will denote by $R_2=\{r\in R\mid 2r=0 \}$.
\begin{trm}\label{FinalTheoremLastOdd}
Let $R$ be a commutative ring. Let $G$ be a
nonabelian group with a nontrivial orientation
homomorphism $\sigma$. Let $N=Ker(\sigma)$ and
denote by $E$ an elementary abelian $2$-group. Then,
$(RG)_{\varphi_{\sigma}}^{-}$ is commutative if and
only if one of the following conditions holds
\begin{enumerate}
\item $R_2=\{0\}$, $G=\langle a,b \mid a^8=1,\; b^2=a^4,\;
ab=ba^3\rangle\times E$ and $N=\langle a^2,ab
\rangle\times E$;
\item   $\Char(R)=4$ and $G=\langle a,b\mid a^8=1,\; b^2=a^4,\;
ab=ba^{-1}\rangle\times E$, $N=\langle a^2,b
\rangle\times E$ or $N= \langle a^2,ab \rangle\times
E$;
\item $R_2=\{0\}$
and $N$ is an elementary abelian $2$-group;
\item  $G$ is a Hamiltonian 2-group and one of the following conditions is satisfied:
\begin{itemize}
\item[(i)] $N$ is abelian,
\item[(ii)] $N$ is a Hamiltonian 2-group and $\Char(R)=4$;
\end{itemize}
\item $G=\GEN{a,b\mid a^4=b^4=1,\; ab=b^{-1}a}\times E$ and
$N$ is equal to either $\GEN{a,b^{2}}$ or
$\GEN{ab,b^2}\times E$,

\item $\Char(R)=4$, $G=\langle a,b,c\ \mid \ a^4=b^4=1,\; c^2=a^2,
\; ab=ba, \; ac=ca^{-1},\;  bc=cb^{-1} \rangle
\times E$ and $N$ is equal to either
$\GEN{a,c}\times \GEN{b^2}\times E$ or
$\GEN{a,bc}\times \GEN{b^2}\times E$;
\item   $R_2=\{0\}$, $G=\langle
a,b,c \mid a^2=b^2=c^2=1,\;  abc=bca=cab
\rangle\times E$ and $N$ is equal to either
$\GEN{a,b}\times E$, $\GEN{a,c}\times E$ or
$\GEN{b,c}\times E$;

\item $R_2=\{0\}$, $G=\langle a,b,c,d \mid
a^4=b^2=c^2=d^2=1,\; ab=ba,\; ac=ca,\; ad=dab,\;
bc=cb,\; bd=db,\; cd=da^2c \rangle\times E$ and
$N=\langle b\rangle\times \langle c,d\rangle\times
E$;
\item $R_2=\{0\}$, $G=\langle a,b,c \mid  a^4=b^4=c^2=1,\; ab=ba,\; ac=ca^{-1},\;
bc=ca^2b^{-1}\rangle\times E$ and $N=\langle
a,c\rangle\times \langle b^2\rangle\times E$;
 \item $R_2=\{0\}$,
$G=\langle a,b,c \mid a^4=b^4=c^2=1,\; ab=ba,\;
ac=ca,\; bc=ca^2b\rangle\times E$ and $N=\langle
b,c\rangle \times E$ or $N=\langle ab,c\rangle
\times E$.
\end{enumerate}
\end{trm}

%
The outline of the paper is as follows. In Section~2
we give several examples, and in particular the
sufficiency of the conditions in the theorem follow.
In Section~3  we prove several technical lemmas. It
follows that if $(RG)_{\varphi_{\sigma}}^{-}$ is
commutative then the exponent of $G$ divides $8$.
In Section~4 we deal with groups of exponent $8$
 (this corresponds with cases 1 and 2 of Theorem~\ref{FinalTheoremLastOdd}). In Section~5 we handle groups
of exponent $4$ and abelian kernel $N$ (this
corresponds with cases 3, 4 with abelian kernel and
5 of Theorem~\ref{FinalTheoremLastOdd}). Finally in
Section~6 the remaining cases are dealt with, that
is, $G$ has exponent $4$ and $N$ is not abelian
(this corresponds with cases 4 with $N$ a
Hamiltonian $2$-group and cases 6 to 10 of
Theorem~\ref{FinalTheoremLastOdd})).

\section{Sufficient Conditions}

In this section we give several examples of finite
groups $G$ with a nontrivial orientation morphism
$\sigma : G\rightarrow \{ -1,1\}$ so that
$(RG)^{-}_{\varphi_{\sigma}}$ is commutative for any
commutative ring $R$. These examples are needed to
prove the sufficiency of the conditions in the main
result.

Throughout $R$ is a commutative ring of
characteristic not $2$ and $G$ is a group with
nontrivial orientation morphism $\sigma$. The
classical involution on $G$ is denoted by $*$. The
order of $g\in G$ is denoted by $\circ (g)$ and the
center of $G$ is denoted by $\z (G)$. For subsets
$X$ and $Y$ of a ring $T$ we denote by $[X,Y]$ the
set of commutators $[x,y]=xy-yx$ with $x\in X$ and
$y\in Y$,
 and the multiplicative commutator
$ghg^{-1}h^{-1}$ of $g,h\in G$ is denoted by
$(g,h)$.

The kernel of $\sigma$ will always be denoted by $N$
and by assumption it always is a proper subgroup of
$G$. So, $N$ is a subgroup of index 2 in $G$. It is
obvious that the involution $\varphi_{\sigma}$
coincides on the subring $RN$ with the ring
involution $*$ and  that the antisymmetric elements
in $G$, under $\varphi_{\sigma}$, are the symmetric
elements in $G\setminus N$ under $*$. Then as an
$R$-module, $(RG)_{\varphi_{\sigma}}^-$ is generated
by the set
\begin{eqnarray*}
{\mathcal S} &=&\{g\in G\setminus
N|\;g^2=1\}\cup\{g-g^{-1} \ | \ g\in N \}
\cup\{g+g^{-1}  \mid
 \ g\in (G\setminus N),\; g^2\neq 1
\}\\ && \;\cup\;\{rg\mid g\in N,\; g^2=1 \mbox{ and
}r\in R_2\}.
\end{eqnarray*}


We begin with stating an obvious but useful remark.

\begin{remark}\label{obsexp}
Let $G=H\times E$, a direct product of groups, with
$E$ an elementary abelian $2$-group. Let $\sigma$
and $\sigma_1$ be orientation homomorphisms of $G$
and $H$ with kernels $N$ and $N_1$, respectively. If
$N=N_1\times E$ then $(RG)_{\varphi_{\sigma}}$ is
commutative if and only if
$(RH)^-_{\varphi_{{\sigma}_1}}$ is commutative.
\end{remark}

Our first example is that of Hamiltonian $2$-groups.

\begin{prop}\label{pca1}
If $G$ is a Hamiltonian $2$-group then
\begin{enumerate}
\item If $N$ is abelian $(RG)_{\vp_\sigma }^{-}$ is commutative.
\item
If $N$ is not abelian and $\Char(R)=4$ abelian then
$(RG)_{\vp_\sigma }^{-}$ is commutative.
\end{enumerate}
\end{prop}
\begin{proof}
1. Assume $N$ is abelian. Because of
Remark~\ref{obsexp}, it is sufficient to deal with
the case $G=Q_8=\GEN{a,b \mid  a^4=1, b^2=a^2,
b^{-1}ab=a^{-1}}$ and $N = \langle a\rangle$.
Because elements of order $2$ are central in $G$, we
only need to check that $[a-a^{-1},x+x^{-1}]=0$ and
$[y+y^{-1},x+x^{-1}]=0$ for all $x,y\not\in N$. For
the former we may assume that $x=b$. Since
$ab=a^{-1}b^{-1}$, $ab^{-1}=a^{-1}b$,
$ba=b^{-1}a^{-1}$ and $ba^{-1}=b^{-1}a$ we get that
\begin{eqnarray*}
[a-a^{-1},b+b^{-1}]&=&ab+ab^{-1}-a^{-1}b-a^{-1}b^{-1}\\
& &-ba+ba^{-1}-b^{-1}a+b^{-1}a^{-1}=0,
\end{eqnarray*}
as desired.

For the remaining case it is sufficient to deal with
$x=b$ and $y=ab$. Now
\begin{eqnarray*}
[ab+b^{-1}a^{-1},b+b^{-1}]&=&ab^2+a+b^{-1}a^{-1}b+b^{-1}a^{-1}b^{-1}\\
& &-bab-a^{-1}-b^{-1}ab-b^{-2}a^{-1}\\
&=&a^{-1}+a+a+a^{-1}-a-a^{-1}-a^{-1}-a=0,
\end{eqnarray*}
again as desired.

2. Assume that $\Char(R)=4$ and $N$ is not abelian,
i.e. it is  Hamiltonian $2$-group. Then $G=N\times
E$ with $E$ a cyclic group of order $2$. It is
easily checked that the antisymmetric elements in
$RN$  commute. This also follows from Example 4.1 in
\cite{BP} (one uses that $\Char (R)=4$). Since $E$
is central in $G$ it then also is easily checked
that $(RG)_{\vp_\sigma }^{-}$ is commutative.
\end{proof}

 Next we deal with four groups of order $16$.
 We will write $G_{[a,b]}$ to denote  the  group  $[a,b]$ in
The Small Group library in GAP \cite{gap}.

\begin{prop}\label{pca2}
Let $G=G_{[16,8]}=\langle a,b \mid  a^8=1,\;
b^2=a^4,\; ab=ba^3\rangle =\langle a\rangle\cup
\langle a\rangle b$ and $R$ a commutative ring with
$R_2=\{0\}$. Then, $N=\langle a^2, ab\rangle=\langle
a^2\rangle\cup\langle a^2\rangle ab$ is the  only
proper kernel for which
$(RG)_{\varphi_{\sigma}}^{-}$ is commutative.
\end{prop}
\begin{proof}
The only subgroups of index $2$ in $G$ are $\langle
a \rangle$, $\langle a^{2},b\rangle$ and $\langle
a^{2},ab\rangle$. In the first two cases we have
that $ab,a^{3}b\in (RG)_{\varphi_{\sigma}}^{-}$, but
these elements do not commute.

So the only possible kernel is  $N=\langle a^2,
ab\rangle=\{1,a^2,a^4,a^6,ab,a^3b,a^5b,a^7b\}$. We
need to show that then $(RG)_{\varphi_{\sigma}}^{-}$
is commutative. Since $R_2=\{0\}$ and the only
elements of order $2$ are  $a^4$, $ab$, $a^3b$,
$a^5b$ and $a^7b$, it is enough to show  that
$[A,A]=0$ with $A=\{a^2-a^6,\; a+a^{-1},\;
a^3+a^5,\; b+b^{-1},\; a^2b+a^6b\}$. Because
$\langle a^{2},b\rangle \cong Q_{8}$ it is
sufficient to check that $[a+a^{-1},b+b^{-1}]=0$. As
$\langle a,b\rangle =\langle a,a^{2}b\rangle$ and
$\circ (a^{2}b)=4$, it follows that
$[a+a^{-1},a^{2}b+(a^{2}b)^{-1}]=[a+a^{-1},a^{2}b+a^{6}b]=0$.
Hence the result follows. Therefore $a\not\in N$.
\end{proof}

\begin{prop}\label{pca3}
 Let
$G=G_{[16,9]}=\langle a,b \mid  a^8=1,\; b^2=a^4,\;
ab=ba^{-1}\rangle =\langle a\rangle\cup \langle
a\rangle b$ and $R$ a commutative ring with
$\Char(R)=4$. Then $\langle a^2, b\rangle$ and
$\langle a^2, ab\rangle$  are the only kernels $N$
for which $(RG)_{\varphi_{\sigma}}^{-}$ is
commutative. Therefore $a\not\in N$.
\end{prop}
\begin{proof}
Let $N=\langle a^2, b\rangle=\langle
a^2\rangle\cup\langle a^2\rangle b$. Since
$\Char(R)=4$ and $N\cong Q_8$, one easily checks
that $(RN)^-$ is commutative (or see \cite{BP}). Let
$A=\{a+a^{-1}, a^3+a^{-3}, ab+ab^{-1},
a^3b+a^{-1}b\}$ and $B=\{ a^2-a^{-2}, b-b^{-1},
a^2b-a^{-2}b\}$. Because  $a^{4}$ is the only
element of order $2$ in $G$ and $a^4$ is central, it
is enough to show that $[A,A]=0$ and $[A,B]=0$.
Clearly the elements that only depend on $a$
commute. Also $[ab+ab^{-1},a^3b+a^{-1}b]=$
$[ab+ab^{-1},a^{-1}b^{-1}+a^{-1}b]=$
$(2-2)a^2+(2-2)a^{-2}=0$.

Now notice that if $g,h\in G$ are such that
$o(g)=8$, $h^2=g^4$ and $gh=hg^{-1}$ then
$[g+g^{-1},h\pm h^{-1}]=[g+g^{-1},h\pm hg^4]=0$.
Thus, since $G=\langle a, ab\rangle=\langle
a,a^3b\rangle= \langle a^3, ab\rangle=\langle a^3,
a^3b\rangle$, we get that $[A,A]=0$. As $G=\langle
a, b\rangle=\langle a,a^2b\rangle= \langle a^3,
b\rangle=\langle a^3, a^2b\rangle$ we obtain that
$0=[a+a^{-1},b-b^{-1}]=$ $[a+a^{-1},a^2b-a^{-2}b]=$
$[a^3+a^{-3},b-b^{-1}]=$
$[a^3+a^{-3},a^2b-a^{-2}b].$ Finally, if $g,h\in G$
are such that $o(g)=4=o(h)$ then $g^2=h^2$.
Therefore $h^{-1}=g^2h=hg^2$ and
$[g+g^{-1},h-h^{-1}]=0$ and hence $[A,B]=0$.

Replacing $b$ by $ab$ we also get the result for
$N=\langle a^2, ab\rangle$.
\end{proof}

\begin{prop}\label{pca4}
 Let
$G=G_{[16,13]}= \langle a,b, c \mid a^2=b^2=c^2=1,\;
abc=bca=cab \rangle= \langle a,b\rangle\cup\langle
a,b\rangle c$ and $R$ a commutative ring with
$R_2=\{0\}$. Then the only kernels $N=\ker (\sigma
)$ for which $(RG)_{\varphi_{\sigma}}^{-}$ is
commutative are $N=\langle a,b\rangle, \GEN{a, c}$
and $\GEN{b,c}$.  (Note that these are all the
nonabelian subgroups isomorphic to $D_{4}$.)
\end{prop}
 \begin{proof}
Note that $\z (G)=\langle abc\rangle =\{
1,abc,(ab)^2,bac\}$ and $G$ is of exponent $4$.

First we show that $(RG)_{\varphi_{\sigma}}^{-}$ is
commutative with $N=\langle a,b \rangle$. Since
$R_2=\{0\}$ and $ab$ and $ba$ are the only elements
of order $4$ in $N$, it is enough to show that
$[ab-ba,A]=0$ and $[A,A]=0$, where $A=\{ c, (ab)^2c,
ac+ca, bc+cb\}$. Let $x=abc\in \z (G)$. Then,
$x^{-1}=bac$ and $ab-ba=(x-x^{-1})c$. Thus clearly
$ab-ba$ commutes with $c$ and $x^2c=(ab)^{2}c$.
Moreover, since $ac+ca=b(bac+bca)=b(x^{-1}+x)$,
$bc+cb=bc+(bc)^{-1}=bc+abac=a(x+x^{-1})$ and
$(x-x^{-1})(x+x^{-1})=x^2-x^{-2}=0$, we have that
$[ab-ba,A]=0$. As $(ab)^2$ is central, it is easy to
see that $c$ and $(ab)^2c$ commute with $ac+ca$ and
$bc+cb$. Finally,
$[ac+ca,bc+cb]=(x+x^{-1})^2[b,a]=(x+x^{-1})^2(x^{-1}-x)c=0$
and $[A,A]=0$. Therefore
$(RG)_{\varphi_{\sigma}}^{-}$ indeed is commutative.

Analogously, due to the symmetry in the generators
of $G$, we have that if $N$ is equal to either
$\langle a,c\rangle$ or $\langle b,c\rangle$ then
$(RG)_{\varphi_{\sigma}}^{-}$ is commutative.

Notice that $G$ has four other possible  kernels:
$N_1=\langle a,bc\rangle$, $N_2=\langle
b,ac\rangle$, $N_3=\langle c,ab\rangle$ and
$N_4=\langle ab,ac\rangle$. If $N=N_1$ then $b,c\in
(RG)_{\varphi_{\sigma}}^{-}$ and they do not
commute. If $N=N_2$ then $a,c\in
(RG)_{\varphi_{\sigma}}^{-}$ and they do not
commute. If $N=N_3$ then $b,a\in
(RG)_{\varphi_{\sigma}}^{-}$ and they do not
commute. Finally if $N=N_4$  then $b, a\in
(RG)_{\varphi_{\sigma}}^{-}$ and they do not
commute.
 \end{proof}

\begin{prop}\label{e16/4}
Let $G=G_{[16,4]}=\langle a,b \mid a^4=b^4=1,\;
ab=b^{-1}a \rangle$. Then, $\langle a, b^2 \rangle$
and $\langle ab,b^2 \rangle$ are the only
 kernels $N=\ker (\sigma )$ for which
$(RG)_{\varphi_{\sigma}}^{-}$ is commutative.
\end{prop}
\begin{proof}
Notice that $\z (G)=\langle a^2 \rangle\times
\langle b^2 \rangle$ and that the only subgroups of
index $2$ in $G$ are $\langle a, b^2 \rangle$,
$\langle ab,b^2 \rangle$ and $\langle a^2,b
\rangle$.

First assume  $N=\langle a, b^2 \rangle=\langle a
\rangle\times\langle b^2 \rangle$. Since $N$ is
abelian and elements of order $2$ are central, to
prove that $(RG)_{\varphi_{\sigma}}^{-}$ is
commutative, it is enough to show that $[A,A]=0$ and
$[A,B]=0$, where  $A=\{ b+b^{-1}, ab+a^{-1}b,
ab^{-1}+a^{-1}b^{-1},a^2b+a^2b^{-1} \}$ and  $B=\{
a-a^{-1}, ab^2-a^{-1}b^2\}$. Now, since
$b+b^{-1}=(1+b^2)b$ and $ab+a^{-1}b=(1+a^2)ab$ we
have that
$[b+b^{-1},ab+a^{-1}b]=(1+a^2)(1+b^2)[b,ab]=(1+
a^2)(1+b^2)a(1-b^2)=$ $a(1+ a^2)(1-b^4)=0$.
Similarly, since
$ab^{-1}+a^{-1}b^{-1}=(1+a^2)ab^{-1}$,
$a^2b+a^2b^{-1}=a^2(1+b^2)b$,
$[b,ab^{-1}]=-a(1-b^2)$ and $[ab,ab^{-1}]=0$ we have
that $[A,A]=0$. On the other hand, since
$a-a^{-1}=a(1-a^2)$, $ab^2-a^{-1}b^2=a(1-a^2)b^2$
and $(1-a^4)=0$ we have that the elements in $[\{
ab+a^{-1}b,\; ab^{-1}+a^{-1}b^{-1} \}, B]=0$.
Moreover, since $[a,b]=ab(1-b^2)$ and $(1-b^4)=0$ we
have that $[\{ b+b^{-1}, \; a^2b+a^2b^{-1}\}, B]=0$.
Thus $[A,B]=0$ and therefore
$(RG)_{\varphi_{\sigma}}^{-}$ indeed is commutative.

Replacing $a$ by $ab$ we also get that
$(RG)_{\varphi_{\sigma}}^{-}$ is commutative for
$N=\langle ab,b^2 \rangle$ then
$(RG)_{\varphi_{\sigma}}^{-}$ is commutative.

On the other hand, if $N=\langle b,a^2
\rangle=\langle b \rangle\times\langle a^2 \rangle$
then $(RG)_{\varphi_{\sigma}}^{-}$ is not
commutative because
\begin{eqnarray*}
[a+a^{-1},b-b^{-1}]&=&(1+a^2)(1-b^2)[a,b]=(1+a^2)(1-b^2)^2ab\\
&=&(1+a^2)2(1-b^2)ab=2(1+ a^2-b^2-a^2b^2)ab\\
&=&2(ab+a^{-1}b-ab^{-1}-a^{-1}b^{-1})\neq 0.
\end{eqnarray*}
\end{proof}

Thirdly we deal with four groups of order $32$.

\begin{prop}\label{e32/35}
Let $G=G_{[32,35]}=\langle a,b,c \mid a^4=b^4=1,\;
c^2=a^2,\; ab=ba,\;  ac=ca^{-1},\; bc=cb^{-1}
\rangle$ and $R$ a commutative ring with
$\Char(R)=4$. Then the only kernels $N=\ker (\sigma
)$ for which $(RG)_{\varphi_{\sigma}}^{-}$ is
commutative are $\langle a, c \rangle\times\langle
b^2 \rangle$ and $\langle a,bc \rangle\times\langle
b^2 \rangle$.
\end{prop}
\begin{proof}
First, notice that $\z (G)=\langle a^2 \rangle\times
\langle b^2 \rangle$. Let $N=\langle a, c
\rangle\times\langle b^2 \rangle$. Then, since
$\Char(R)=4$ and $N\cong Q_8\times C_2$, we have
that $(RN)_{\varphi_{\sigma}}^{-}$ is commutative
(see \cite{BP}). Let $A_0=\{b+b^3\}$,
$A_1=\{a^3b+ab^3, bc+a^2bc, abc-a^3bc\}$,
$B_0=\{a-a^3, c-a^2c, ac+a^3c\}$, $A=A_0\cup
A_0a^2\cup A_1\cup A_1b^2$ and $B=B_0\cup B_0b^2$.
To prove that $(RG)_{\varphi_{\sigma}}^{-}$ is
commutative, it is enough to show
$[A,A]=[B,B]=[A,B]=0$.  Since $a^2$ and $b^2$ are
central, it is enough to show $[A_{0}\cup
A_{1},A_{0}\cup A_{1}]=[B_{0},B_{0}]=[A_{0}\cup
A_{1},B_{0}]=0$. Again, as $a^2$ and $b^2$ are
central, we can write $b+b^3=(1+b^2)b$,
$a^3b+ab^3=(a^2+b^2)ab,$ $bc+a^2bc=(1+a^2)bc,$ and
$abc+a^3bc=(1+ a^2)abc$. Thus, since $ab=ba$, it
follows that $[b+b^3,a^3b+ab^3]=0$. Also,
$[b+b^3,bc+a^2bc]=(1+a^2)(1+b^2)[b,bc]=(1+a^2)(1+b^2)(b^2-1)c=0$.
Similarly, as  $[b,abc]=(b^2-1)ac$,
$[ab,bc]=(b^2-a^2)ac$ and
$[ab,abc]=a^{2}b^{2}c-abcab=(a^{2}b^{2}-1)c$, we get
that $[b+b^{3},abc+a^3bc]=0$ and
 $[a^3b+ab^3,\{ bc+a^2bc,\; abc+a^3bc\}]=0$.
Moreover, as $[bc,abc]=(1-a^2)a$ and
$(1+a^2)(1-a^2)=0$, we obtain  that
$[bc+a^2bc,abc+a^3bc]=0.$ Therefore $[A_0\cup A_1 ,
A_{0}\cup A_{1}]=0 $.

Because $\Char (R)=4$, we also obtain that
$(1-a^{2})^{3}=0$. It then easily follows that the
elements of $[B_{0},B_{0}]=0$.

On the other hand, since $a-a^3=(1-a^2)a$, we have
that $[a-a^3,A_0\cup A_1]=0$. Also, as
$c-a^2c=(1-a^2)c$, $[c,b]=(b^2-1)bc$,
$[c,ab]=(a^2b^2-1)abc$ and $(1+b^2)(1-b^2)=0=(1+
a^2b^2)(1-a^2b^2)$, it follows that $[c-a^2c,
A_0\cup A_1 ]=0$. Similarly, as $ac-a^3c=(1-a^2)ac,$
$[ac,b]=(b^2-1)abc$ and $[ac,ab]=(b^2-a^2)bc$, we
have that $[ac-a^3c, A_0\cup A_1 ]=0$. Therefore the
elements of $[A_{0}\cup A_{1},B_0]=0$.  Hence
$(RG)_{\varphi_{\sigma}}^{-}$ is commutative.

Similarly, replacing $c$ by $bc$, we get that if
$\Char(R)=4$ and $N=\langle a,bc\rangle\times\langle
b^2\rangle$ then $(RG)_{\varphi_{\sigma}}^{-}$ is
commutative.

Notice that $G$ has five other possible kernels:
$N_1=\langle a,b\rangle$, $N_2=\langle a^2,b,c
\rangle$, $N_3=\langle a^2,b,ac \rangle$,
$N_4=\langle a^2, a^3b,c \rangle$ and $N_5=\langle
a^2, a^3b, ac \rangle$. But, if $N=N_1$, $N_2$ or
$N_3$ then $b-b^3\in (RG)_{\varphi_{\sigma}}^{-}$
and either $c-a^2c\in (RG)_{\varphi_{\sigma}}^{-}$
or $c+a^2c\in (RG)_{\varphi_{\sigma}}^{-}$, and
$[b-b^3,c\mp a^2c]=-2(1\mp a^2)(1-b^2)bc\neq 0$. On
the other hand, if $N=N_4$ or $N_5$ then
$ab-a^3b^3\in (RG)_{\varphi_{\sigma}}^{-}$ and
either $c-a^2c\in (RG)_{\varphi_{\sigma}}^{-}$ or
$c+a^2c\in (RG)_{\varphi_{\sigma}}^{-}$, and
$[ab-a^3b^3,c\mp a^2c]=2(1\mp a^2)(1-a^2b^2)abc\neq
0$.
\end{proof}

\begin{prop}\label{e32/38}
Let $G=G_{[32,30]}=\langle a,b,c,d \mid
a^4=b^2=c^2=d^2=1,\; ab=ba,\; ac=ca,\; ad=dab,\;
bc=cb,\; bd=db,\; cd=da^2c \rangle$ and $R$ a
commutative ring with $R_2=\{0\}$. Then, $N=\langle
b \rangle\times\langle c,d \rangle$ is the only
kernel $N$ for which $(RG)_{\varphi_{\sigma}}^{-}$
is commutative.
\end{prop}
\begin{proof}
Notice $\z (G)=\langle a^2 \rangle\times \langle b
\rangle$. Assume $N=\langle b \rangle\times\langle
c,d \rangle$. Let $A_0=\{a+a^3, ac+a^3c,
ad+a^3bd\}$, $A_1=\{acd+abcd\}$, $B_0=\{cd-a^2cd\}$,
$A=A_0\cup A_0 b\cup A_1\cup A_1 a^2$ and $B=B_0\cup
B_0b$. Since $R_2=\{0\}$, to prove that
$(RG)_{\varphi_{\sigma}}^{-}$ is commutative, it is
enough to show that $[A,A]=[B,B]=[A,B]=0$. Clearly,
as $b$ is central, the elements of $[B,B]=0$.
Moreover, since $a^2$ is central, we can write
$a+a^3=(1+a^2)a$, $ac+a^3c=(1+a^2)ac$ and
$acd+abcd=(1+ b)acd$. Thus, since $ac=ca$, we get
that $[a+a^3,ac+a^3c]=0$. Also, as
$[a,acd]=(1-b)a^2cd$ and $(1+b)(1-b)=0$, we have
that $[a+a^3,acd+acdb]=0$. Now,
$[a+a^3,ad+a^3bd]=(1+a^2)(1+b)[a,ad]=(1+a^{2})(1+b)(1-b)a^2d=0$.
Similarly, as $(1+a^2)(1+a^2b)=(1+b)(1+a^2b)$,
$[ac,ad]=(a^2-b)cd$, $[ac,acd]=(a^2-b)d$ and
$[ad,acd]=(1-a^2)bc$, we have that
$[a+a^3,ad+a^3bd]=0$, $[ac+a^3c, \{ad+a^3bd,\;
acd+abcd]=0$ and $[ad+ a^3bd,acd+abcd]=0$. Thus,
$[A_0\cup A_1, A_{0}\cup A_{1}]=0$ and thus
$[A,A]=0$. Since $cd-a^2cd=(1-a^2)cd$,
$0=(1+a^2)(1-a^2)=(1+a^2b)(1-a^2b)=(1-b)(1+b)$,
$[ad,cd]=a^3c(1-a^2b)$ and $[acd,cd]=a^3(1-b)$, we
have that $[A_{0}\cup A_{1},B_{0}]=0$. Therefore the
elements of $[A,B]=0$. Hence
$(RG)_{\varphi_{\sigma}}^{-}$ is commutative.

Now, notice that $G$ has six other possible kernels:
$N_1=\langle a,b,c \rangle$, $N_2=\langle
a^2,b,c,a^3bd \rangle$, $N_3=\langle a,b,d \rangle$,
$N_4=\langle a^2, b, a^3c,d \rangle$, $N_5=\langle
a, b, a^2cd \rangle$ and $N_6=\langle a^2, b, a^3c,
a^3bd \rangle$. If $N=N_1$ or $N_2$ then $d\in
(RG)_{\varphi_{\sigma}}^{-}$ and either $a-a^3\in
(RG)_{\varphi_{\sigma}}^{-}$ or $a+a^3\in
(RG)_{\varphi_{\sigma}}^{-}$, but $[a\mp
a^3,d]=(1\mp a^2)(b-1)da\neq 0$. If $N=N_3$ or $N_4$
then $c\in (RG)_{\varphi_{\sigma}}^{-}$ and either
$ad-a^3bd\in (RG)_{\varphi_{\sigma}}^{-}$ or
$ad+a^3bd\in (RG)_{\varphi_{\sigma}}^{-}$, but
$[c,ad\mp a^3bd]=(1\mp a^2b)(1-a^2)acd\neq 0$.
Finally, if $N=N_5$ or $N_6$ then $c,d\in
(RG)_{\varphi_{\sigma}}^{-}$ and they do not
commute.
\end{proof}

\begin{prop}\label{e32/39}
Let $G=G_{[32,31]}=\langle a,b,c \mid
a^4=b^4=c^2=1,\; ab=ba,\; ac=ca^{-1},\;
bc=ca^2b^{-1}\rangle$ and $R$ a commutative ring
with $R_2=\{0\}$. Then, $\langle a,c
\rangle\times\langle b^2 \rangle$ is the only
 kernel  for which $(RG)_{\varphi_{\sigma}}^{-} $ is
commutative.
\end{prop}
\begin{proof}
Notice that $\z (G)=\langle a^2 \rangle\times
\langle b^2 \rangle$. Assume $N=\langle a,c
\rangle\times\langle b^2 \rangle$. Let
$A_0=\{b+b^3\}$, $A_1=\{ab+a^3b^3,\; bc+a^2bc,\;
abc+a^3bc\}$, $B_0=\{a-a^3\}$, $A=A_0\cup A_0a^2\cup
A_1\cup A_1b^2$ and $B=B_0\cup B_0b^2$. Since
$R_2=\{0\}$, to prove that
$(RG)_{\varphi_{\sigma}}^{-}$ is commutative it is
enough to show that $[A,A]=[B,B]=[A,B]=0$. Clearly,
as $b^2$ is central, it follows that $[B,B]=0$.
Since $a^2$ is  central, we can write $A_{1}=\{
ab(1+a^2b^2),\; (1+a^2)bc, \; (1+a^2)abc\}$. Thus,
because $ab=ba$, it follows that
$[b+b^3,ab+a^3b^3]=0$. Also, as $[bc,abc]=a(1-a^2)$
and $(1+a^2)(1-a^2)=0$, we have that
$[bc+a^2bc,abc+a^3bc]=0$. On the other hand, as
$(1+a^2)(1+b^2)=(1+a^2)(1+a^2b^2)$,
$[b,bc]=(b^2-a^2)c$, $[b,abc]=(b^2-a^2)ac$,
$[bc,ab]=0$ and $[ab,abc]=(b^2-1)a^2c$, we have that
$[A_0\cup A_1 , A_{0}\cap A_{1}]=0$ and thus
$[A,A]=0$.

As $a-a^3=(1-a^2)a$, $ab=ba$ and $(1+a^2)(1-a^2)=0$,
we have $[A_{0}\cup A_{a},B_{0}]=0$. Therefore
$[A,B]=0$  and hence $(RG)_{\varphi_{\sigma}}^{-}$
is commutative.

Notice that $G$ has six other possible  kernels:
$N_1=\langle a^2,b,c \rangle$, $N_2=\langle
a^2,a^2b,c \rangle$, $N_3=\langle a^2,b,ac \rangle$,
$N_4=\langle a^2, a^3b, ac \rangle$, $N_5=\langle a,
b^2, b^2c \rangle$ and $N_6=\langle a, b\rangle$. If
$N=N_1$ or $N_2$ then $ac\in
(RG)_{\varphi_{\sigma}}^{-}$ and either $b-b^3\in
(RG)_{\varphi_{\sigma}}^{-}$ or $b+b^3\in
(RG)_{\varphi_{\sigma}}^{-}$, but $[ac,b\mp
b^3]=(1\mp b^2)(a^2b^2-1)abc\neq 0$. If $N=N_3$ or
$N_4$ then $c\in (RG)_{\varphi_{\sigma}}^{-}$ and
and either $b-b^3\in (RG)_{\varphi_{\sigma}}^{-}$ or
$b+b^3\in (RG)_{\varphi_{\sigma}}^{-}$, but $[c,b\mp
b^3]=(1\mp b^2)(a^2b^2-1)bc\neq 0$. If $N=N_5$ or
$N_6$ then $c,ac\in (RG)_{\varphi_{\sigma}}^{-}$ and
they do not commute.
\end{proof}

\begin{prop}\label{e32/16}
Let $G=G_{[32,24]}=\langle a,b,c \mid
a^4=b^4=c^2=1,\; ab=ba,\; ac=ca,\; bc=ca^2b\rangle$
and $R$ a commutative ring with $R_2=\{0\}$. Then,
the only kernels  for which
$(RG)_{\varphi_{\sigma}}^{-}$ is commutative are
$\langle b,c \rangle$ and $\langle ab,c \rangle$.
\end{prop}
\begin{proof}
Notice that $\z (G)=\langle a \rangle\times \langle
b^2 \rangle$. Assume $N=\langle b,c \rangle$. Then,
since $N$ contains an elementary abelian
$2$-subgroup of index $2$, it follows that $(RN)^-$
is commutative (see \cite{BP}). Now, let
$A_0=\{ac+a^3c, ab+a^3b^3\}$, $A_1=\{abc+ab^3c\}$,
$B_0=\{b-b^3\}$, $B_1=\{bc-a^2b^3c\}$, $A=A_0\cup
A_0b^2\cup A_1\cup A_1a^2$ and $B=B_0\cup B_0a^2\cup
B_1\cup B_1b^2$. Since $R_2=\{0\}$, to prove that
$(RG)_{\varphi_{\sigma}}^{-}$ is commutative, it is
enough to show $[A,A]=[A,B]=[B_{0},B_{1}]=0$. The
last equality follows from
$[(1-b^{2})b,(1-a^{2}b^{2})bc]=
(1-b^{2})(1-a^{2}b^{2})[b,bc] =
(1-b^{2})(1-a^{2}b^{2}) (a^{2}-1) c = 0$. Since
$a^2$ and $b^2$ are central, we can write
$ab+a^3b^3=a(1+ a^2b^2)b$, $ac+a^3c=a(1+ a^2)c$,
$abc+ab^3c=a(1+ b^2)bc$, $b-b^3=(1-b^2)b$ and
$bc-a^2b^3c=(1- a^2b^2)bc$. Thus $[ab+a^3b^3,
b-b^3]=0$ and, as $(1+a^2b^2)(1-a^2b^2)=0$,
$[ab+a^3b^3,bc-a^2b^3c]=0$. Also, since
$[b,bc]=(1-a^2)b^2c$ and $(1+a^2b^2)(1+
b^2)(1-a^2)=0$, we have that $[ab+a^3b^3 =
(1+a^{2}b^{2})ab,abc+ab^3c = (1+b^{2})abc] = 0$. On
the other hand, since $(1+a^2b^2)(1-a^2b^2)=0$,
$[c,b]=(a^2-1)bc$ and $[c,bc]=(a^2-1)b$, we have
that $[ac+a^3c, \{ ab+a^3b^3,\; abc+ab^3c,\;
b-b^3,\; bc-a^2b^3c\}]=0$. Finally, notice that
$[abc+ab^3c, bc-a^2b^3c]=0$ and, as $(1+
b^2)(1-b^2)=0$, $[abc+ab^3c,b-b^3]=0$. Therefore,
$[A_0\cup A_1, A_{0}\cup A_{1}]=0$ and $[A_0\cup
A_1, B_0\cup B_1]=0$. Hence
$(RG)_{\varphi_{\sigma}}^{-}$ is commutative.

Notice that $G=\langle a,ab,c\rangle$, $o(ab)=4$,
$a(ab)=(ab)a$ and $(ab)c=ca^2(ab)^{-1}$. Hence,
replacing $b$ by $ab$, we get that also $N=\langle
ab,c\rangle$ is a  kernel.

Finally,  $G$ has five other possible kernels:
$\langle a,b^2,c \rangle$, $\langle a,b \rangle$,
$\langle a^2,b,a^3c \rangle$, $\langle a, b^2,
a^2b^3c \rangle$ and $\langle a^2, a^3b, a^3c
\rangle$. If $N=\langle a,b^2,c\rangle$ then
$ac-a^3c$, $b+b^3\in (RG)_{\varphi_{\sigma}}^{-}$
and they do not commute. Otherwise, $c\in
(RG)_{\varphi_{\sigma}}^{-}$ and either $b-b^3\in
(RG)_{\varphi_{\sigma}}^{-}$ or $b+b^3\in
(RG)_{\varphi_{\sigma}}^{-}$, but $[c, b\pm
b^3]=(a^2-1)(1\pm b^2)bc\neq 0$.
\end{proof}

%
%

We finish this section with one more elementary
remark.

\begin{remark}\label{cenpro}
Let $G$ be a group and let $A$ be a subgroup of
index 2 in $G$. Assume that $A=C\times E$, a direct
product of groups, with $E$  an elementary abelian
$2$-group. If $E$ is central in $G$ then $G$ is the
central product of the subgroups $E$ and $\langle C,
g\rangle$, with $g\in G\setminus A$.\end{remark}


 \section{Necessary conditions}\label{I}

We begin with a series of technical lemmas that
yield necessary conditions for
$(RG)_{\varphi_{\sigma}}^{-}$ to be commutative.
\begin{lem}\label{fixedcentral}
Assume that $(RG)_{\varphi_{\sigma}}^{-}$ is
commutative.
\begin{enumerate}
\item If $R_2=\{0\}$ and  $g\in G\setminus N$ with $g^2=1$ then $(g,h)=1$ for all
$h\in G\setminus N$ with $h^2=1$ and for all $h\in
N$ with $h^2\neq 1$.
\item If $R_2\neq\{0\}$  and  $g\in G$ with $g^2=1$ then $g\in\z(G)$.
\end{enumerate}
\end{lem}
\begin{proof}
1.  Assume that $R_2=\{0\}$ and  $g\in G\setminus N$
with $g^2=1$.
 If $h\in G\setminus N$ with $h^2=1$, then $g$ and $h$ are two antisymmetric elements and hence $[g,h]=0$ as desired.
Assume now that $h\in N$ with $h^2\neq 1$. Then
$0=[g,h-h^{-1}]=gh-gh^{-1}-hg+h^{-1}g$ and therefore
 we have that $gh$ is equal to either $gh^{-1}$ or $hg$. The former is excluded
since by assumption $h^{2}\neq 1$. Thus $gh=hg$, as
desired.

2. Assume that $g\in G$ with $g^2=1$ and $0\neq r\in
R_2$. Suppose first that $g\in N$. Let $h\in G$. We
need to show that $(g,h)=1$. If $h\in G\setminus N$
with $h^2=1$ then $0=[rg,h]=rgh-rhg$ and therefore
$gh=hg$. If $h\in G\setminus N$ with $h^2\neq 1$
then $0=[rg,h+h^{-1}]=r(gh+gh^{-1}+hg+h^{-1}g)$.
Since $h^2\neq 1$, we clearly have that $gh\neq
gh^{-1}$. Hence either $gh=hg$ (as desired) or
$gh=h^{-1}g$. The latter implies that
 $(gh)^2=1$ with $gh\in G\setminus N$. So,
by the previous case, $1=(g, gh)=(g,h)$ as desired.
Finally if  $h\in N$ then choose $x\in G\setminus
N$. By the previous, $(g,x)=1=(g,hx)$. Hence
$1=(g,h)$ as desired.

Second, assume that $g\in G\setminus N$. Let $h\in
G\setminus N$. If $h^2=1$ then, by the above,
$0=[g,h]$, as desired. If $h^2\neq 1$ then
$0=[g,h+h^{-1}]=gh+gh^{-1}-hg-h^{-1}g$. So either
$gh=hg$ (as desired) or $gh=h^{-1}g$. The latter
implies that$(gh)^2=1$ and $gh\in N$. So, by the
above, $1=(gh,h)$ and thus $1=(g,h)$. We thus have
shown that $g$ commutes with all elements $h\in
G\setminus N$. Assume now that $h\in N$. If $h^2=1$
then by the above $1=(g,h)$, as desired. If $h^2\neq
1$ then $0=[g, h+h^{-1}]=gh+gh^{-1}-hg-h^{-1}g$. It
follows that either $gh=hg$ (as desired) or
$gh=h^{-1}g$ and thus $(gh)^2=1$ with $gh\in
G\setminus N$ and by the above $1=(g,gh)=(g,h)$
which finishes the proof of the lemma.
\end{proof}

\begin{lem}\label{ca0}
Let $g$ and $h$ be elements of $G$ with $g^2\neq 1$
and  $h^2\neq 1$.  If $[g-\varphi_{\sigma}(g),
h-\varphi_{\sigma}(h)]=0$ then the following
properties hold.
\begin{enumerate}
\item[$(i)$]
If $g,h\in N$ then one of the following conditions
holds.
 \begin{enumerate}
  \item $gh=hg$.
  \item $R_2=\{0\}$ and $(g^\alpha h^\beta)^2=1$, for all $\alpha , \beta \in \{ -  1,1\}$.
  \item $\Char(R)=4$ and $\langle g,h\rangle\cong Q_8$.
 \end{enumerate}
\item[$(ii)$]
If $g\in N$ and $h\not\in N$ then one of the
following conditions
              holds.
 \begin{enumerate}
  \item $ghg^{-1}\in \{ h,h^{-1}\}$.
  \item $\circ(g)=4=\circ(h)$ and $g^2=h^2$.
 \end{enumerate}
\item[$(iii)$] If $g,h\not\in N$ then one of the following conditions holds:
 \begin{enumerate}
  \item $gh\in \{ hg,h^{-1}g,hg^{-1}\}$.
  \item $R_2=\{0\}$ and $(g^\alpha h^\beta)^2=1$, for all $\alpha , \beta \in
        \{ -1,1\}$.
 \end{enumerate}
\end{enumerate}
\end{lem}

\begin{proof} $(i)$ By \cite[Lemma 2.1]{BP}  we have that either $gh=hg$; or $(g^\alpha h^\beta)^2=1$, for all
 $\alpha , \beta \in \{ -  1,1\}$; or $\Char(R)=4$ and $\langle g,h\rangle\cong Q_8$. Notice that
if $(g^\alpha h^\beta)^2=1$, for all
 $\alpha , \beta \in \{ -  1,1\}$ and $R_2\neq \{0\}$, then by Lemma~\ref{fixedcentral} it follows that
  $1=(gh,h)=(g,h)$ so we  are in case $(a)$.

$(ii)$ Suppose $g\in N$ and $h\in G\setminus N$.
Then,  $0=[g-g^{-1},h+h^{-1}]=
gh+gh^{-1}-g^{-1}h-g^{-1}h^{-1}-hg+hg^{-1}-h^{-1}g+h^{-1}g^{-1}.$
As $g^2\neq 1$,   $h^2\neq 1$  and $\Char(R)\neq 2$,
it follows that $gh$ equals either $hg$, $h^{-1}g$,
or $g^{-1}h^{-1}$.

Assume that $gh=g^{-1}h^{-1}$, that is,
$g^2=h^{-2}$. Then
$0=gh^{-1}-g^{-1}h+hg^{-1}-h^{-1}g$ and so $gh^{-1}$
is equal to either $g^{-1}h$ or $h^{-1}g$.
Therefore, $g^2=h^2$ or $gh=hg$. But, if $g^2=h^2$
then we obtain that $\circ(g) = 4 = \circ(h)$. Hence
(ii) follows.

$(iii)$ Suppose $g,h\not\in N$. Then,
\begin{equation}\label{as}
0=[g+g^{-1},h+h^{-1}]=
gh+gh^{-1}+g^{-1}h+g^{-1}h^{-1}-hg-hg^{-1}-h^{-1}g-h^{-1}g^{-1}.
\end{equation}
As $g^2\neq 1$,  $h^2\neq 1$  and $\Char(R)\neq 2$,
it follows that $gh$ equals  either $hg$, $h^{-1}g$,
$hg^{-1}$, or $h^{-1}g^{-1}$.

Assume that $gh=h^{-1}g^{-1}$, that is, $(gh)^2=1$,
or equivalently $(g^{-1}h^{-1})^2=1$. If
$R_2\neq\{0\}$ then, by Lemma~\ref{fixedcentral}, it
follows that $(g,h)=1$.
 So assume that $R_2=\{0\}$. By (\ref{as}) we know that
$0=gh^{-1}+g^{-1}h-hg^{-1}-h^{-1}g$. Thus, either
$gh^{-1}=hg^{-1}$ or $gh^{-1}=h^{-1}g$. Therefore,
$(gh^{-1})^2=1$ and $(g^{-1}h)^2=1$, or $gh=hg$.
This finishes the proof of the lemma.
\end{proof}

If $(RG)_{\varphi_{\sigma}}^{-}$ is commutative then
the following remark can be applied to elements of
order $2$ that do not belong to $N$.

%

\begin{lem}\label{ca1}
Assume that $(RG)_{\varphi_{\sigma}}^{-}$ is
commutative with $\Char(R)\neq 2$. Let $g$ and $h$
be noncommuting elements of $G$ such that $g^2\neq
1$ and $ h^2\neq 1$. The following properties hold.
\begin{enumerate}
\item[$(i)$] If $g,h\in N$ then one of
      the following conditions holds.
 \begin{enumerate}
  \item $R_2=\{0\}$ and $\langle g,h\rangle =  \langle g,h\ |\
 g^4=h^4=(gh)^2=(gh^{-1})^2=1\rangle=G_{[16,3]}$.
  \item $\Char(R)=4$ and $\langle g,h\rangle\cong Q_8$.
 \end{enumerate}
\item[$(ii)$] If $g\in N$ and $h\not\in N$ then one of the following conditions
              holds.
 \begin{enumerate}
   \item $\langle g,h\rangle\cong G_{[16,4]}$.
   \item $\langle g,h\rangle\cong Q_8$;
   \item $\Char(R)= 4$ and  $\langle g,h\rangle \cong G_{[16,9]}$.
 \end{enumerate}
\item[$(iii)$] If $g,h\not\in N$ then one of the following conditions
holds.
 \begin{enumerate}
\item $\langle g,h\rangle$ is isomorphic to either $Q_8$,  or $G_{[16,4]}$.
  \item $R_2=\{0\}$ and $\langle g,h\rangle$ is isomorphic to either $G_{[16,3]}$, or $G_{[16,8]}$.
  \item $\Char(R)= 4$ and  $\langle g,h\rangle\cong G_{[16,9]}$.
 \end{enumerate}
\end{enumerate}
\end{lem}

\begin{proof}  Assume $g,h\in G$ are noncommuting and
$g^{2}\neq 1$, $h^{2}\neq 1$.

$(i)$ Suppose $g,h\in N$. Because of Lemma~\ref{ca0}
$(i)$ we may assume that $R_2=\{0\}$ and
$(gh)^2=(gh^{-1})^2=1$. Hence, it remains to show
that $\circ(g)= 4$ and $\circ(h) = 4$. We prove the
former, the latter is similar.

Let $x\in G\setminus N$. We may assume that $\circ
(x)\neq 2$. Indeed, for assume $\circ (x)=2$ then,
by Lemma~\ref{fixedcentral} $(1)$, $gx=xg$. Thus
$(gx)^2=g^2\neq 1$ and we may replace $x$ by $gx \in
G\setminus N$.

First we deal with the case that $gx\neq xg$.
Applying Lemma \ref{ca0} $(ii)$ to the elements $g$
and $x$, we get that either $gx=x^{-1}g$, or $\circ
(g)=4$ and $g^2=x^2$. If $gx=x^{-1}g$  then,
applying Lemma \ref{ca0} $(ii)$ to the elements $g$
and $gx$, it follows that either
$g^2x=(gx)^{-1}g=x^{-1}$ and hence $g^2=x^{-2}$, or
$\circ (g)=4=\circ (gx)$ and $g^2=(gx)^2$.
The former implies that
$x^{-2}=gx^{2}g^{-1}=g^{-2}=x^{2}$ and thus
$g^{4}=x^{-4}=1$, i.e. $\circ (g)=4$. So we have
shown that if $gx\neq xg$ then $\circ (g)=4$.

Second we deal with the case that $gx=xg$. Because,
by assumption $hg\neq gh$ and thus $g(hx)\neq
(hx)g$, Lemma~\ref{fixedcentral} $(1)$ yields that
$(hx)^2\neq 1$. Thus, Lemma \ref{ca0} $(ii)$ applied
to $g$ and $hx$, gives that either $ghx=hxg=hgx$,
$ghx=(hx)^{-1}g$, or $\circ (g)=4$. The former is
excluded because $gh\neq hg$. Since,
$(gh)^{2}=(gh^{-1})^{2}=1$ we also know that
$ghx=h^{-1}g^{-1}x=h^{-1}xg^{-1}$ and
$x^{-1}h^{-1}g=x^{-1}g^{-1}h=g^{-1}x^{-1}h$, the
second option thus implies that $ghx=x^{-1}h^{-1}g=
h^{-1}xg^{-1}$ and also
$ghx=x^{-1}h^{-1}g=g^{-1}x^{-1}h$. Hence
$g^{-2}=x^{-1}hx^{-1}h^{-1}$ and
$g^{2}=x^{-1}hx^{-1}h^{-1}$. So $g^{4}=1$.
Therefore, we have shown that $gx=xg$ implies that
$\circ (g)=4$. This finishes the proof of (i).

$(ii)$ Suppose $g\in N$ and $h\not\in N$.  Because
of  Lemma~\ref{ca0} $(ii)$, we know that either
$gh=h^{-1}g$, or $\circ(g)=4=\circ(h)$ and
$g^2=h^2$.

First, suppose that $gh=h^{-1}g$ and so $g^2$ is a
central element in the group $\langle g,h\rangle$.
Since $(gh)^2=g^{2}\neq 1$ we can  apply
Lemma~\ref{ca0} $(ii)$ to the elements $g$ and $gh$,
and we obtain that either
$g^{2}h=(gh)^{-1}g=h^{-1}$, or $\circ (g)=4$. Thus,
either $g^2=h^{-2}$ or $\circ (g)=4$. On the other
hand, applying Lemma~\ref{ca0} $(i)$ to $g$ and
$h^{-1}gh$, we have that either
$gh^{-1}gh=h^{-1}ghg$; or  $R_2=\{0\}$ and
$(g^{-1}h^{-1}gh)^2=1$, or $\Char(R)=4$ and $\langle
g,h^{-1}gh\rangle\cong Q_{8}$. Consequently, since
$gh=h^{-1}g$, we get that either
$g^2h^2=h^{-2}g^{2}=g^2h^{-2}$; or $R_2=\{0\}$ and
$h^{4}=1$; or $\Char(R)=4$ and
$g^{2}=(gh^{-1}gh)^{2}=(g^{2}h^{2})^{2}=g^{4}h^{4}$.
Therefore, either $h^4=1$, or $\Char(R)=4$ and
$g^2=h^{-4}$. Hence $\langle g,h\rangle$ is
isomorphic to either $Q_8$, $G_{[16,4]}$, or
$\Char(R)=4$ and $\langle g,h\rangle\cong
G_{[16,9]}$.

Second, suppose that $\circ(g)=4=\circ(h)$ and
$g^2=h^2$. Lemma~\ref{ca0} $(i)$ applied to $g$ and
$hgh$ yields that  either $(gh)^2=(hg)^2$; or
$R_2=\{0\}$ and $(gh)^4=1$; or  $\Char(R)=4$ and
$\langle g,hgh\rangle\cong Q_{8}$. Thus, either
$(gh)^4=1$ or $\Char(R)=4$ and $g^2=(gh)^4$.
Therefore, either $\langle g,h\rangle \cong Q_{8}$,
$\langle g,h\rangle= \langle g,h\ \mid
g^4=(gh)^4=1,\; g^2=h^2 \rangle\cong  G_{[16,4]}$,
or $\Char(R)=4$ and $\langle g,h\rangle = \langle
g,h\ \mid \ g^4=1,\; (gh)^4=g^2=h^2 \rangle \cong
G_{[16,9]}$.


$(iii)$ Suppose $g,h\not\in N$.   By Lemma~\ref{ca0}
$(iii)$ we have that either $gh=h^{-1}g$,
$gh=hg^{-1}$, or  $(g^\alpha h^\beta)^2=1$ for all
$\alpha,\beta \in \{ -1,1\}$ and $R_2=\{0\}$.

First, suppose that $gh=h^{-1}g$ and so $g^2$ is
central in $\langle g,h\rangle$. Then, applying
Lemma~\ref{ca0} $(ii)$ to $gh$ and $g$, we obtain
that either $ghg=g^{-1}gh=h$ and hence $g^2=h^2$, or
$\circ (g)=4$. On the other hand, applying
Lemma~\ref{ca0} $(i)$ to $gh$ and $hg$, we get that
either $gh^2g=hg^2h$; or $R_2=\{0\}$ and
$(gh(hg)^{-1})^2=1$; or $\Char(R)=4$ and $\langle
gh,hg\rangle\cong Q_8$. Thus, we have that either
$g^2h^{-2}=g^2h^2$; or $R_2=\{0\}$ and $h^4=1$; or
$\Char(R)=4$, $1=(gh)^4=g^4$ and
$(gh)^2=(hg^2h)^2=(g^2h^2)^2=g^4h^4$. Therefore,
either $h^4=1$, or $\Char(R)=4$, $g^2=h^4$ and
$g^{4}=1$. Hence, either $\langle g,h\rangle \cong
Q_8$, $\langle g,h\rangle\cong G_{[16,4]}$, or
$\Char(R)=4$ and $\langle g,h\rangle\cong
G_{[16,9]}$.

Second, suppose $gh=hg^{-1}$. Then the result
follows at once from the previous case by replacing
$g$ by $gh$ and $h$ by $g^{-1}$.

Third, suppose that $R_2=\{0\}$ and $(g^\alpha
h^\beta)^2=1$, for all $\alpha,\beta \in \{ -
1,1\}$. In particular, $hg^{-1}h=g$ and
$ghg=h^{-1}$.  Lemma~\ref{ca0} $(iii)$ applied  to
$ghg^{-1}$ and $h$ yields that $((ghg^{-1})^\alpha
h^\beta)^2=1$ for all $\alpha,\beta \in \{ -1,1\}$
(and, in particular, $(ghg^{-1}h)^2=(g^2)^2=g^4=1$)
or $ghg^{-1}h=g^2$ is equal to either
$hghg^{-1}=g^{-2}$, $h^{-1}ghg^{-1}=h^{-2}g^{-2}$ or
$hgh^{-1}g^{-1}=h^2g^{-2}$. Thus, either $\circ
(g)=4$ or $g^4=h^{\pm 2}$. Similarly, applying
Lemma~\ref{ca0} $(iii)$ to $g$ and $hgh^{-1}$, it
follows that either $\circ (h)=4$ or $h^4=g^{\pm
2}$.

If $\circ (g)=4$ and $\circ (h)=4$ then $\langle
g,h\rangle \cong G_{[16,3]}$. On the other hand, if
$\circ (g)=4$ and $h^4=g^2$, or $\circ (h)=4$ and
$g^4=h^2$, then $\langle g,h\rangle \cong
G_{[16,8]}$. Now, assume that  $\circ (g)\neq 4$ and
$\circ (h)\neq 4$.  Then $g^4=h^{\pm 2}$ and
$h^4=g^{\pm 2}$ and hence $g^2=g^4g^{-2}=h^{\pm
2}h^{\mp 4}=h^{\mp 2}$. So $g^2 h^{\pm 2}=1$. On the
other hand, since $(g^\alpha h^\beta)^2=1$, for all
$\alpha,\beta \in \{ - 1,1\}$, we have that
$1=g^2h^{\pm 2}=gh^{\mp 1}g^{-1}h^{\pm 1}$.
Therefore $gh^{\mp 1}=h^{\mp 1}g$, and thus $gh=hg$,
a contradiction.
\end{proof}

\begin{lem}\label{ca2}
Assume that $(RG)_{\varphi_{\sigma}}^{-}$ is
commutative. Let $g$ and $h$ be elements of $G$ with
$g^2\neq 1$ and $h^2=1$. If  $R_2=\{0\}$ then the
following properties hold.
\begin{enumerate}
\item[$(i)$] If $g,h\in N$ then one of the following  conditions
holds.
  \begin{enumerate}
    \item $\langle g,h\rangle$ is abelian.
    \item  $(gh)^2\neq 1$ and $\langle g,h\rangle \cong
           G_{[16,3]}$.
    \item  $(gh)^2=1$ and $\langle g,h\rangle=
           \langle g,h\ |\ g^4=h^2=(gh)^2=1\rangle=D_4$.
  \end{enumerate}
\item[$(ii)$]
If $g\in N$ and $h\not\in N$ then $\langle
g,h\rangle$ is abelian.
\item[$(iii)$] If $h\in N$ and $g\not\in N$ then one of the following conditions
holds.
  \begin{enumerate}
   \item $\langle g,h\rangle$ is abelian.
   \item $(gh)^2\neq 1$ and $\langle g,h\rangle\ $ is isomorphic to either
         $G_{[16,3]}$ or $G_{[16,8]}$.
   \item $(gh)^2=1$and $\langle g,h\rangle\cong D_4$.
  \end{enumerate}
\item[$(iv)$] If $g,h\not\in N$ then  $\langle
g,h\rangle$either  is abelian or isomorphic to
$D_4$.
\end{enumerate}
If $R_2\neq \{0\}$ then $\GEN{g,h}$ is abelian.
\end{lem}

\begin{proof}
Note that the last part of the statement follows at
once from Lemma~\ref{fixedcentral}. So we assume
throughout the proof that $R_2=\{0\}$. $(i)$ Suppose
$g,h\in N$. Assume that $gh\neq hg$. If $(gh)^2\neq
1$ then, by Lemma \ref{ca1} $(i)$ it follows that
$\langle g,h\rangle=\langle g,gh\rangle\cong
G_{[16,3]}$.

So, to prove $(i)$ we assume from now on that
$(gh)^{2}=1$ and thus $gh=hg^{-1}$.  Choose $x\in
G\setminus N$. We may assume that $x^{2}\neq 1$.
Indeed, for otherwise, by Lemma~\ref{fixedcentral},
$gx=xg$ and thus $(gx)^{2}=g^{2}x^{2}=g^{2}\neq 1$;
so we can replace $x$ by $gx$. We now claim that
$\circ (g)=4$ and therefore $\langle g,h\rangle\cong
D_4$, as desired. We prove this by contradiction.
Hence, assume $\circ (g)\neq 4$. Lemma~\ref{ca0}
$(ii)$ applied to the elements $g\in N$ and $x\in
G\setminus N$ yields that either $gx=xg$ or
$gx=x^{-1}g$.

First, assume that $gx=xg$. Then $(hx)^2\neq 1$,
because otherwise, by Lemma~\ref{fixedcentral}
$(1)$, it follows that $ghx=hxg=hgx$, and hence
$gh=hg$, a contradiction. Thus, because by
assumption $\circ (g)\neq 4$,  applying
Lemma~\ref{ca0} $(ii)$, to $g$ and $hx$, we get that
either $ghx=hxg=hgx$  or $ghx=(hx)^{-1}g$. The
former is excluded as $gh\neq hg$. So
$ghx=x^{-1}hg$. Hence, since $ghx=hg^{-1}x=hxg^{-1}$
and $x^{-1}hg=x^{-1}g^{-1}h=g^{-1}x^{-1}h$, we
obtain that $(hx)^2=g^2=(x^{-1}h)^2$ and therefore
$g^4=(hx)^2(x^{-1}h)^2=1$. This gives a
contradiction with the assumption $\circ (g)\neq 4$.

Second, assume that  $gx=x^{-1}g$. Since $g^2\neq 1$
we have that $(x^{-1}gx)^2\neq 1$. So, applying
Lemma~\ref{ca1} $(i)$ to the elements $g$ and
$x^{-1}gx$ (recall that $\circ (g)\neq 4$), we get
that  $g$ and $x^{-1}gx$ commute. So
$gx^{-1}gx=x^{-1}gxg$ and thus $g^2x^{2}=x^{-2}g^2$.
Now, if $(hx)^2=1$ then, by Lemma~\ref{fixedcentral}
$(1)$, it follows that $ghx=hxg$. Thus
$hg^{-1}x=hgx^{-1}$ and hence $g^2=x^2$. Then
$g^2x^{2}=x^{-2}g^2=1$. Hence $g^2=x^2=g^{-2}$ and
therefore $\circ (g)=4$, a contradiction. So
$(hx)^2\neq 1$ and we can apply Lemma~\ref{ca0}
$(ii)$ to $g$ and $hx$. It follows that either
$ghx=hxg$ or $ghx=(hx)^{-1}g$. We already have shown
above that the former leads to a contradiction.
Hence,  $ghx=x^{-1}hg$.  Since
$ghx=hg^{-1}x=hx^{-1}g^{-1}$ and
$x^{-1}hg=x^{-1}g^{-1}h=g^{-1}xh$, this yields that
$xh(hx)^{-1}=g^2=(x^{-1}h)^{-1}hx^{-1}$ and so
$g^4=(x^{-1}h)^{-1}hx^{-1}xh(hx)^{-1}=1$. Hence
$\circ (g)=4$, again a contradiction. This finishes
the proof of $(i)$.

$(ii)$ This follows at once from
Lemma~\ref{fixedcentral} $(1)$.

$(iii)$ Suppose that $h\in N $, $g\not\in N$ and
$gh\neq hg$. First assume that  $(gh)^2\neq 1$. Then
we can apply Lemma~\ref{ca0} $(iii)$ to $g$ and
$gh$. It follows that either
$(g^\alpha(gh)^\beta)^2=1$ for all $\alpha,\beta\in
\{  -1,1\}$,  $ggh=g^2h=ghg^{-1}$ or
$ggh=(gh)^{-1}g=h$. The latter is excluded as it
yields  $g^{2}=1$. The second possibility leads to
$(gh)^{2}=1$ and is thus also excluded. It follows
from Lemma~\ref{ca1} $(iii)$ that  $\langle
g,h\rangle=\langle g,gh\rangle$ is isomorphic to
either $G_{[16,3]}$ or $G_{[16,8]}$, as desired.

Second assume that $(gh)^2=1$. We claim that then
$\circ (g)=4$, and thus $\langle g,h\rangle\cong
D_4$. Indeed, suppose the contrary, that is $g^4\neq
1$. We then can apply part $(ii)$ to $g^2$ and $gh$
and we get  $g^2gh=ghg^2=gg^{-2}h$. Thus $g^4=1$, a
contradiction.

$(iv)$ Suppose  $g,h\not\in N$. If $(gh)^2\neq 1$
then part $(ii)$ yields that $\langle
g,h\rangle=\langle h,gh\rangle$ is abelian. On the
other hand, if $(gh)^2=1$ then  part $(iii)$ implies
that $\langle g,h\rangle =\langle g,gh\rangle$
either is abelian or $\langle g,h \rangle\cong D_4$,
because $(ggh)^{2}=1$
\end{proof}

\begin{lem}\label{ca3}
Assume that $(RG)_{\varphi_{\sigma}}^{-}$ is
commutative. Let $g$ and $h$ be non-commuting
elements of $G$. If $g^2 = h^2 =1$ then  $R_2=\{0\}$
and $\langle g,h \rangle \cong D_{4}$.
\end{lem}

\begin{proof} Assume that $gh\neq hg$. Hence, as $g^{2}=h^{2}=1$, we get
that $(gh)^{2}\neq 1$. Since also
$(g(gh))^{2}=h^{2}=1$ and $\langle g,h\rangle
=\langle g,gh\rangle$, the result follows from
Lemma~\ref{ca2}.
\end{proof}

\begin{remark} \label{remarkcentral}  Lemmas~\ref{ca1}, \ref{ca2} and
\ref{ca3} imply that if $G$ is a group of exponent
$4$ and $(RG)_{\varphi_{\sigma}}^{-}$ is commutative
for a nontrivial orientation morphism then $g^2\in
\z (G)$, for all $g\in G$. Thus $G'\subseteq \z
(G)$.
\end{remark}

We end this section by showing that if $G$ is a
nonabelian group with  $(RG)_{\varphi_{\sigma}}^{-}$
commutative then $G$ is a $2$-group of  exponent at
most $8$.

\begin{prop}\label{propo}
Let $R$ be a commutative ring with $\Char(R) \neq 2$
and let $G$ be nonabelian group. If
$(RG)_{\varphi_{\sigma}}^{-}$ is commutative then
$G$ is a $2$-group and its exponent is bounded by 8.
\end{prop}

\begin{proof}
It follows from  Lemmas \ref{ca1}, \ref{ca2} and
\ref{ca3} that noncentral elements $g$ of $G$ have
order a divisor of $8$. If $y$ is a central element
of $G$ then $yg$ is a noncentral element and thus
$1=(yg)^{8}=y^{8}g^{8}=y^{8}$. Hence $o(y)$ divides
$8$ as well and the result follows.
\end{proof}

\section{Groups of exponent eight}

We know from Proposition~\ref{propo} that if
$(RG)_{\varphi_{\sigma}}^{-}$ is commutative then
$G$ is a $2$-group of exponent bounded by $8$. In
this chapter we give a complete answer in case the
exponent is precisely $8$.

\begin{trm}\label{exp8}
Let $R$ be a commutative ring with $\Char(R) \neq 2$
and let $G$ be a nonabelian group of exponent $8$
with a nontrivial orientation homomorphism $\sigma$.
Then, $(RG)_{\varphi_{\sigma}}^{-}$ is commutative
if and only if one of the following conditions
holds.
\begin{enumerate}
\item[$(i)$]  $R_2=\{0\}$, $G= \langle g,h \mid g^8=1,\;
h^2=g^4,\; gh=hg^3\rangle\times E$ and $N=\langle
g^2,gh \rangle\times E$, for some elementary abelian
$2$-group  $E$.
\item[$(ii)$] $\Char(R)=4$, $G=\langle
g,h\mid g^8=1,\; h^2=g^4,\; gh=hg^{-1}\rangle\times
E$ and $N=\langle g^2,h \rangle\times E$ or $N=
\langle g^2,gh \rangle\times E$, for some elementary
abelian $2$-group  $E$.
\end{enumerate}
\end{trm}

\begin{proof} Suppose $(RG)_{\varphi_{\sigma}}^{-}$ is commutative and
$G$ is a $2$-group of exponent $8$. Let $A=\langle
a\in G \mid \circ(a)=8\rangle$. By assumption, $A$
is nontrivial.  Fix $a\in A$ with $\circ (a) =8$.
Because of Lemma~\ref{ca1}, we know that $A$ is an
abelian group. We claim that the elements of order
$8$ of $G$ belong to $G\setminus N$, and that
$h^{-1}gh=g^{3}$ or $g^{-1}$, for all $g\in A$ and
$h\in G\setminus A$.

Since $A$ is an abelian group  generated by elements
of order $8$, it is enough to prove the claim for
$g=a$ and $h\in G\setminus A$.  Since $\circ (h)\leq
4$ and $\circ (a)=8$,  note that $ah\neq ha$.
Indeed, because otherwise $\circ (ah)=8$ while
$ah\not \in A$, a contradiction. Lemma~\ref{ca1} and
Lemma~\ref{ca2} then yield that $\langle a,h\rangle
=G_{[16,8]}$ or $G_{[16,9]}$, and by
Propositions~\ref{pca2} and \ref{pca3} also that
$a\not\in N$. In particular we obtain that
$h^{-1}ah=a^{-1}$ or $a^{3}$, as desired. This
finishes the proof of the claim.

Now, we show that $A$ has index $2$ in $G$. In order
to show this, let $x, y\in G\setminus A$. Suppose
that $xy^{-1}\in G\setminus A$. Then, by the
previous paragraph, $a^{-2} =
(xy^{-1})^{-1}a^{2}(xy^{-1}) = y(x^{-1}a^{2}x)y^{-1}
= ya^{-2}y^{-1} = a^{2}$. Hence,  $a^{4}=1$, a
contradiction. Therefore, $xy^{-1}\in A$, and thus
indeed $[G:A]=2$.

Next we show that $A$ is the direct product of a
cyclic group of order $8$ and an elementary abelian
$2$-group. For this, first recall that every abelian
group of finite exponent is a direct product of
cyclic groups of prime power order (see for example
\cite[(5.1.2), p.92]{scot}). Because $A$ is abelian
of exponent $8$, we thus get that  $A$ has a cyclic
subgroup of order $8$ as a direct factor. Without
loss of generality, we may assume that this factor
is $\langle a \rangle$.
In order to show that $A$ does not have a direct
factor that is a cyclic group of order $4$, it is
sufficient to prove that $a^{4}\in \{ c^{2},c^{4}\}$
for any element $c \in A$ with $c^2\neq 1$. So, let
$c\in A$ with $c^{2}\neq 1$. Suppose that $a^4\neq
c^4$. Then $(ac)^4=a^4c^4\neq 1$ and $\circ (ac)=8$.
So, by the claim above, $ac\in G\setminus N$. As
$a\not\in N$ and $[G:N]=2$, we therefore obtain that
$c\in N$. So, again by the above claim, $\circ
(c)=4$. Now, as in the beginning of the proof, let
$h\in G\setminus A$. Then,  $ah\neq ha$. As $G$ is a
$2$-group and $\circ (h)\leq 4$, we may assume that
$\circ (h)=4$. Indeed, for otherwise, $h^2=1$ and by
Lemma~\ref{fixedcentral} it follows that
$R_2=\{0\}$. Then by Lemma~\ref{ca2} we have that
$h\in N$, $\circ (ah)=4$ and $ah\in G\setminus A$.
So, replacing $h$ by $ah$ we obtain the desired.
Then, by Lemma~\ref{ca1}, it follows that
$R_2=\{0\}$ and $\langle a,h \rangle$ is isomorphic
to  $G_{[16,8]}$, or $\Char(R)=4$ and $\langle a,h
\rangle\cong G_{[16,9]}$. Consequently, $a^4=h^2$.
On the other hand, from the claim in the beginning
of the proof we also know that $h^{-1}ch=c^{-1}$.
Part $(i)$ and $(ii)$ of Lemma~\ref{ca0} then yield
that either $ch=hc=c^{-1}h$, or
$c^{-1}h=h^{-1}c=c^{-1}h^{-1}$, or
$ch=h^{-1}c=c^{-1}h^{-1}$, or $c^2=h^2$. Since
$c^2\neq 1$ and $ h^2\neq 1$, we deduce that
$c^2=h^2=a^4$, as desired. So, $A=\langle a
\rangle\times E$, with $E^2=1$.

Notice that by the first part of the proof
$h^{-1}eh=e$ for all $e\in E$. Hence $E$ is central
in $G$. Hence, from Remark~\ref{cenpro}, $G$ is the
central product of $\langle a,h \rangle$ and $E$.
Moreover, from the previous, either $R_2=\{0\}$ and
$\langle a,h \rangle\cong G_{[16,8]}$ or
$\Char(R)=4$ and $\langle a,h \rangle\cong
G_{[16,9]}$. Furthermore, as $\langle a\rangle\cap
E=\{1\}$, either $ah=ha^3$ or $ah=ha^{-1}$, and
$hA=G\setminus A$, we have that $\langle a,
h\rangle\cap E=\{1\}$. Hence $G=\langle
a,h\rangle\times E$.

To finish the proof of the necessity of the
conditions, it remains to determine the kernels. By
Remark~\ref{obsexp} and Propositions \ref{pca2} and
\ref{pca3} we get the desired kernels and also the
sufficiency of the conditions follows.
\end{proof}

\section{Groups of exponent four and abelian kernel}

In the remainder of the paper we are left to deal
with nonabelian $2$-groups $G$ of exponent $4$. In
this section we handle  such groups for which the
kernel $N$ is abelian. Without specific reference to
Remark~\ref{remarkcentral} we will often use he fact
that $g^2\in Z(G)$ for $g\in G$ if
$(RG)^{-}_{\varphi_{\sigma}}$ is commutative.

We first prove that if $N$ is an elementary abelian
$2$-group then $(RG)^{-}_{\varphi_{\sigma}}$ is
commutative.

\begin{prop}\label{N2ele}
Let $G$ be a nonabelian group of exponent $4$, $R$ a
commutative ring with $R_2=\{0\}$ and $\sigma$ a
nontrivial orientation homomorphism.  Assume that
$N$ is an elementary abelian $2$-group. Then
$(RG)_{\varphi_{\sigma}}^{-}$ is commutative.
\end{prop}
\begin{proof} Since $N$ is of index $2$ and elementary
abelian, the nonabelian group $G$ contains an
element $x$ so that $G=N\cup xN$ and $x$ has order
$4$. Furthermore, since $R_2=\{0\}$, to prove the
result, it is sufficient to show that
$[g-\varphi_{\sigma}(g), h-\varphi_{\sigma}(h)]=0$
for all $g,h\in G\setminus N$. Write $g=xa$ and
$h=xb$ for some $a,b\in N$. Then
$$\begin{array}{lcl}
  [g-\varphi_{\sigma} (g), h-\varphi_{\sigma} (h)] &=&
  [g+g^{-1},h+h^{-1}]   \\
   & = &[xa+ax^{-1},xb+bx^{-1}] \\
  & = & xaxb+xabx^{-1}+ab+ax^{-1}bx^{-1} \\
  & & -xbxa-xbax^{-1}-ba-bx^{-1}ax^{-1} \\
  &=& xaxb+ax^{-1}bx^{-1}+-xbxa-bx^{-1}ax^{-1}
\end{array}$$
Let $a',b'\in N$ so that $ax=xa'$ and $bx=xb'$.
Since $G$ has exponent $4$, $x^{2}\in N$ and $N$ is
abelian, we get that
$$[g+g^{-1},h+h^{-1}]=x^2a'b+ab'x^{-2}-x^2b'a-ba'x^{-2}=0,$$
as desired.
\end{proof}

Next, assume $N$ is abelian but not an elementary
abelian $2$-group. The following lemma deals with
elements of order $2$ in $N$.

\begin{lem}\label{Nabel}
Let $R$ be a commutative ring of $\Char(R)\neq 2$,
let $G$ be a nonabelian group of exponent $4$ and
$\sigma$ a nontrivial orientation homomorphism.
Assume that $(RG)^{-}_{\varphi_{\sigma}}$ is
commutative and $N$ is abelian but not elementary
abelian $2$-group. Let $a\in N$. Then, $a^2=1$ if
and only if $a\in \z (G)$. Furthermore, $G_{[16,3]}$
is not a subgroup of $G$ and if $x\in G\setminus N$
then $x$ has order $4$.
\end{lem}

\begin{proof}  Assume that
$(RG)^{-}_{\varphi_{\sigma}}$ is commutative and $N$
is abelian but not elementary abelian $2$-group.

First, we show that if  $x\in G\setminus N$ then
$\circ (x)=4$. Assume the contrary, that is assume
$x\in G\setminus N$ and $x^{2}=1$. Then, by
Lemma~\ref{ca2}, we have that $ax=xa$, for all $a\in
N$ with $a^2\neq 1$. Because of the assumptions, $N$
is generated by elements of order $4$. Hence, we get
that $x$ is central and thus $G$ is abelian, a
contradiction.

Second, we show that  $G_{[16,3]}$ is not a subgroup
of $G$. Assume, the contrary. That is, suppose, that
$H =G_{[16,3]}= \langle g,h\ |\
g^4=h^4=(gh)^2=(gh^{-1})^2=1\rangle \subseteq G$.
Clearly, $N\cap H$ is an abelian subgroup of index
$2$ in $H$. Since $\circ (gh) =\circ (gh^{-1})=2$,
we know from the above that $gh\in N\cap H$. As $N$
is abelian,  $(g,h)\neq 1$ and thus $(g,gh)\neq 1$,
$(h,gh)\neq 1$, we thus get that  $g\not\in N$ and
$h\not\in N$.

Take $a\in N$ with $a^2\neq 1$. Because of
Lemma~\ref{ca0} $(ii)$, we get that  either $ag=ga$,
$ag=g^{-1}a$ or $a^2=g^2$,  and either $ah=ha$,
$ah=h^{-1}a$ or $a^2=h^2$. But $g^2\neq h^2$. Also,
since $gh\in N$, we have that $ag=ga$ is equivalent
to $ah=ha$. Thus $ag=ga$ implies $ah\neq h^{-1}a$
and $a^2\neq h^2$. Indeed, for otherwise we obtain
that $h^2=1$ or $(ah)^2=1$. The former obviously is
false. Because of the first part of the proof, the
latter implies that $ah\in N$, again a
contradiction.  Similarly, $ah=ha$ implies that
$ag\neq g^{-1}a$ and $a^2\neq g^2$. So, we have to
consider four remaining cases: $ag=ga$ and $ah=ha$,
$ag=g^{-1}a$ and $ah=h^{-1}a$, $ag=g^{-1}a$ and
$a^2=h^2$, or $a^2=g^2$ and $ah=h^{-1}a$. We show
that each case leads to a contradiction.

Case 1: $ag=ga$ and $ah=ha$. Then, applying
Lemma~\ref{ca0} $(iii)$ to the elements $g$ and
$ah$, we have that $g(ah)=agh=ah^{-1}g^{-1}$ is
equal to either $(ah)g$, $(ah)^{-1}g$, $(ah)g^{-1}$
or $(ah)^{-1}g^{-1}$. Thus either $gh=hg$, $g^2=1$,
$h^2=1$ or $a^2=1$, a contradiction.

Case 2:  $ag=g^{-1}a$ and $ah=h^{-1}a$. Since
$gh^{-1}\in N$, we then have that
$agh^{-1}=gh^{-1}a=gah=ag^{-1}h$. Hence
$gh^{-1}=g^{-1}h$ and therefore $g^2=h^2$, a
contradiction.

Case 3: $ag=g^{-1}a$ and $a^2=h^2$. Then, since
$gh\in N$, we have that $gha=agh=g^{-1}ah$ and thus
$g^2=g^{-2}=hah^{-1}a^{-1}=(ha)^2$. Hence, applying
Lemma~\ref{ca0} $(iii)$ to the elements $g$ and $ah$
we have that $g(ah)=ag^{-1}h=ah^{-1}g$ is equal to
either $(ah)g$, $(ah)^{-1}g$, $(ah)g^{-1}$ or
$(ah)^{-1}g^{-1}$. So, either $h^2=1$, or
$ah^{-1}=h^{-1}a^{-1}=ha$ and then
$g^2=haha=a^2=h^2$, or  $h^{-1}g=hg^{-1}=gh^{-1}$
and hence $gh=hg$, or
$g^2=ha^{-1}h^{-1}a^{-1}=(ha)^2a^2=g^2a^2$ and then
$a^2=1$. Therefore, each of the possibilities yields
a contradiction.

Case 4:  $a^2=g^2$ and $ah=h^{-1}a$. Similarly as in
Case 3, applying Lemma~\ref{ca0} $(iii)$ to the
elements $ag$ and $h$, we obtain a contradiction.

So, indeed we have shown that $G$ does not have
$G_{[16,3]}$ as a subgroup.

Now, assume $a\in N$ with $\circ(a)=2$. If $a\not\in
\z (G)$ then  by Lemma~\ref{fixedcentral} it follows
that $R_2=\{0\}$. Moreover, there
 exists $x\in G\setminus N$  such that
$ax\neq xa$ and $x^2\neq 1$. As $G$ has exponent $4$
and $G_{[16,3]}\not\subset G$, Lemma \ref{ca2}
$(iii)$ yields that $\langle a,x\rangle\cong D_4$.
Then $\circ(ax)=2$ and thus, by the first part of
the proof $ax\in N$. This of course is impossible.
Hence, we have shown that elements of order $2$ in
$N$  are central in $G$. It remains to show that the
converse holds. We prove this by contradiction. So,
assume  $a\in \z (G)\cap N$ and  $a^2\neq 1$. Since
$N$ is an abelian subgroup of exponent $4$ and of
index $2$ in $G$, and because $G$ is not abelian,
there exists $b\in G\setminus N$ and $c\in N$ such
that $bc\neq cb$ and $\circ (c)=4$. Again by the
first part of the proof, $\circ(b)=4$. Note also
that $\circ (ba)=4$.  By Lemma~\ref{ca0} $(ii)$, we
have that $cb=b^{-1}c$ or $c^2=b^2$. Assume first
that $cb=b^{-1}c$. Then, applying Lemma \ref{ca0}
$(ii)$ to the elements $ba$ and $c$, we get that
$b^2a^2=c^2$. Hence $a^2\neq c^2$, as $b^{2}\neq 1$.
Consequently, $(ac)^2\neq 1$. Thus, applying
Lemma~\ref{ca0} $(ii)$ to the elements $ba$ and
$ca$, we obtain that $caba=b^{-1}a^{-1}ca=b^{-1}c$
or $b^2a^2=c^2a^2$. Because $a^{2}\neq 1$, the
former is excluded. However, because  $b^2a^2=c^2$,
the latter also yields that $a^{2}=1$, again a
contradiction. Hence,  $b^2=c^2$. Notice that, if
$a^2=c^2=b^2$ then $(ba)^2=1$ and, by
Lemma~\ref{ca2} $(ii)$, we thus have that $bac=cab$
and therefore $bc=cb$, a contradiction. So,
$(ac)^2\neq 1$. Lemma~\ref{ca0} $(ii)$, applied  to
the elements $b$ and $ac$, then yields that
$acb=b^{-1}ac$ or $b^2a^2=c^2$. However the latter
is impossible as  $a^2\neq 1$ and $b^{2}=c^{2}$.
Therefore $c^2=b^2$ and $cb=b^{-1}c$. Now, applying
Lemma~\ref{ca0} $(ii)$ to the elements $ba$ and $c$,
we obtain that $cba=a^{-1}b^{-1}c$ or $b^2a^2=c^2$.
Both cases imply that $a^2=1$,  a contradiction.
This finishes the proof of the Lemma.
\end{proof}

We are  now in a position the prove a solution to
the problem in case the kernel is abelian.

\begin{trm}\label{exp4Nabelian}
Let $R$ be a commutative ring with $\Char(R) \neq 2$
and let $G$ be a nonabelian group of exponent $4$
with a nontrivial orientation homomorphism $\sigma$.
Assume that $N$ is abelian. Then,
$(RG)^-_{\varphi_{\sigma}}$ is commutative if and
only if one of the following conditions holds.
\begin{enumerate}
\item[$(i)$] $R_2=\{0\}$ and $N$ is an elementary abelian $2$-group.
\item[$(ii)$] $G\cong Q_8\times E$ and $N=C_4\times E$, where $C_4$ is a cyclic group of order 4 and $E$ is
an  elementary abelian $2$-group.
\item[$(iii)$]
$G=\GEN{a,b\mid a^4=b^4=1, ab=b^{-1}a}\times E$ and
$N=\GEN{a,b^{2}} \times E$ or
$N=\GEN{ab,b^{2}}\times E$,
where $E$ is an elementary abelian $2$-group.
\end{enumerate}
\end{trm}
\begin{proof}
Assume $(RG)_{\varphi_{\sigma}}^{-}$ is commutative.
Since $G$ is not abelian, if $N$ is elementary
abelian 2-group then $R_2= \{0\}$,  by
Lemma~\ref{fixedcentral}. So, suppose that $N$ is
not an elementary abelian $2$-group. We need to show
that either (ii) or (iii) holds.

First, suppose that $G_{[16,4]}\not\subseteq G$.
Then, by Lemma~\ref{fixedcentral}, Lemma~\ref{ca1},
Lemma~\ref{ca2} and Lemma~\ref{Nabel}, we have that
$G$ is a Hamiltonian $2$-group, that is, $G\cong
Q_8\times E$, where $E^2=1$. Because $N$ has index
$2$, it is then also clear that $N=C_{4}\times E$
for some elementary abelian $2$-subgroup $E$ of $G$.

Second, suppose that $G_{[16,4]}=\langle g,h \ |\
g^4=h^4=1, gh=h^{-1}g\rangle\subseteq G$. Then,
since $N$ is  abelian, we have that $g\in G\setminus
N$ or $h\in G\setminus N$. We claim that $h\in
G\setminus N$. Indeed, for suppose the contrary.
Then $h\in N$. Hence, by Lemma~\ref{ca0} $(ii)$, we
have that $hg=g^{-1}h$. Since $gh=h^{-1}g$, one
deduces that $g^2=h^2$, a contradiction. This
finishes the proof of the claim.

Put $a=g$ if $g\not\in N$, otherwise put $a=gh$.
Clearly, $a\not\in N$, $a^{2}=g^{2}$, $ah=h^{-1}a$,
$\circ(a)=4$ and $\GEN{a,h}=\GEN{g,h}$. So $N$ is an
abelian group of exponent $4$ and it contains $ah$.
We claim that $N=\GEN{ah}\times E$ for some
elementary abelian $2$-group $E$.  For this it is
sufficient to show that  if $c\in N$ with $c^{2}\neq
1$ then $\GEN{ah}\cap \GEN{c}\neq \{ 1\}$. Suppose
the contrary. Then $ah$, $c$ and $ahc$ have order
$4$. Hence, because of Lemma~\ref{Nabel}, $ah$,  $c$
and $ahc$ are not central in $G$. Because $N$ has
index $2$ in $G$ and since $h\not\in N$, we get that
$ch\neq hc$ and $ahch\neq hahc$. Hence, applying
Lemma~\ref{ca0} $(ii)$  to the elements $ac$ and
$h$, we get that $ahch=h^{-1}ahc$ or
$(ahc)^{2}=h^{2}$. Because  $ah=h^{-1}a$, the former
is excluded. Hence,
$h^{2}=(ahc)^{2}=(ah)^{2}c^{2}=a^2c^2$. Applying
Lemma~\ref{ca0} $(ii)$ to the elements $c$ and $h$,
we also obtain that $ch=h^{-1}c$ or $c^{2}=h^{2}$.
As $a^{2}\neq 1$ and $h^{2}=a^{2}c^{2}$, it follows
that $ch=h^{-1}c$ must hold.  So,
$ahch=ahh^{-1}c=ac=hh^{-1}ac=hahc$ and thus $ahc\in
\z (G)$. Since also $ahc\in N$, Lemma~\ref{Nabel}
yields that $(ahc)^2=1$,  in contradiction with
$\circ (ahc)=4$. This finishes the proof of the
claim that $N=\GEN{ah}\times E$ for some elementary
abelian $2$-group $E$. Again by Lemma~\ref{Nabel},
we also know that $E\subseteq \z (G)$.

Because $(ah)^{2}=a^{2}=g^{2} \neq h^{2}$, it is
clear that $N=\GEN{ah}\times \GEN{h^{2}}\times E_0$,
for some elementary abelian subgroup $E_0$ of $E$.
Note that $\GEN{ah,h^{2}}$ equals either
$\GEN{g,h^{2}}$ or $\GEN{gh,h^{2}}$. Moreover, since
the only central elements of order $2$ in
$\GEN{a,h}=\GEN{g,h}$ are $g^{2}$, $h^{2}$ and
$g^{2}h^{2}$ and since none of these belong to $E$,
we also get that $G=\GEN{g,h}\times E_0$. This
finishes the proof of the necessity of the
conditions.

The sufficiency of the conditions follows from
Remark~\ref{obsexp}, Proposition~\ref{pca1},
Proposition~\ref{e16/4} and Proposition~\ref{N2ele}.
\end{proof}


\section{Groups of exponent four and nonabelian
kernel}

In this section we handle the remaining case, that
is,  we consider groups $G$ of exponent four and
with nonabelian kernel $N$. We first solve our
problem in case all elements of order $2$ in $N$ are
central in $N$.

\begin{lem}\label{ord2cen}
Let $R$ be a commutative ring with $\Char(R) \neq 2$
and let $G$ be a nonabelian group of exponent  $4$.
Assume that $N$  is  not abelian and that the
elements of order 2 in $N$ are central in $N$. If
$(RG)^{-}_{\varphi_{\sigma}}$ is commutative and
$x\in G$ with $x^2=1$ then $x\in \z (G)$.
Furthermore, $\Char(R)=4$, $N$ is a Hamiltonian
$2$-group and $G_{[16,3]}$ is not a subgroup of $G$.
\end{lem}
\begin{proof} Assume $(RG)^{-}_{\varphi_{\sigma}}$ is
commutative. If $R_2\neq \{0\}$ then by
Lemma~\ref{fixedcentral} the first part of the
result follows. So assume that  $R_2=\{0\}$ and fix
$h$ and $h_1$ in $N$ so that $(h,h_{1})\neq 1$.
Because of the assumptions, $\circ (h)=\circ
(h_{1})=4$. Let $x\in G$ with $\circ (x)=2$. We need
to show that $x\in \z (G)$, or equivalently,
$(x,g)=1$ for all $g\in G$.

First, assume  $x\not\in N$. Then, by
Lemma~\ref{ca2} $(ii)$, for all $g\in N$ with
$g^2\neq 1$,  we have that $gx=xg$. Note that, in
particular, $hx=xh$. Now, consider $g\in N$ with
$g^2=1$. Because, by assumption, $g$ is central in
$N$, we get that $(gh)^{2}=g^{2}h^{2} = h^{2}\neq
1$. Hence, $gh$ commutes with $x$ and thus
$xgh=ghx=gxh$. So, $xg=gx$. We thus have shown that
$(x,g)=1$ for all $g\in N$. Since $N$ has index $2$
in $G$ and $x\not\in N$,  we get that $x\in \z (G)$.

Second, assume that $x\in N$. By assumption, $x\in
\z (N)$. Let $g\in G\setminus N$. In order to prove
that $x\in \z (G)$, it is sufficient to show that
$gx=xg$. If $g^{2}=1$ then this follows from the
previous. If $(gx)^{2}=1$ or $(gh)^{2}=1$ then,
again by the above, $gx, \, gh\in \z (G)$. In the
former case, $gx=xg$. In the latter case,
$xgh=ghx=gxh$ and thus $xg=gx$, as desired. So, we
may assume that $\circ (g)=\circ (gx) =\circ
(gh)=4$.  By Lemma~\ref{ca0} $(ii)$, for $y\in N$
with $\circ (y)=4$, we then have three
possibilities: (1) $gy=yg$, (2) $g^2=y^2$ or (3)
$yg=g^{-1}y$. Of course, this can be applied to the
elements $y=h$ or $y=h_{1}$. It is therefore
sufficient to consider the following three cases.

(1) $gh=hg$ (or, by symmetry, $gh_{1}=h_{1}g$).
Lemma~\ref{ca0} $(ii)$, applied to the elements $g$
and $hx$, yields that either $hxg=ghx=hgx$,
$hxg=g^{-1}hx=hg^{-1}x$ or $g^{2}=(hx)^{2}=h^{2}$.
So, either $xg=gx$, $(gx)^{2}=gxgx=xg^{-1}gx=1$ or
$(gh)^{2}=g^{2}h^{2}=1$. Because the latter two are
excluded, we get that $xg=gx$, as desired.

(2) $g^2=h^2$ and $gh\neq hg$ (or, by symmetry,
$g^{2}=h_{1}^{2}$ and $gh_{1}\neq h_{1}g$). Because
$xh=hx$, we have that $gxh\neq hgx$ and therefore,
by Lemma~\ref{ca0} $(ii)$, we get that either
$hgx=xg^{-1}h$ or $(gx)^2=h^2$. If $(gx)^2=h^2=g^2$
then $gx=xg$. Therefore, we may assume that
$hgx=xg^{-1}h$. Lemma~\ref{ca0} $(ii)$, applied the
elements $gx$ and $hx$,  yields that either
$hxgx=gxhx$ or $hxgx=xg^{-1}hx$ or
$(gx)^2=(hx)^2=h^2$. In the first case we have that
$hxg=gxh=ghx$. Hence,
$gh=hxgx=xhgx=xxg^{-1}h=g^{-1}h$ and therefore
$g=g^{-1}$, a contradiction. In the second case, we
obtain $hxgx=xg^{-1}hx=hgxx=hg$ and hence $gx=xg$,
as desired. In the third case, we have
$(gx)^2=h^2=g^2$ and thus  $gx=xg$, again as
desired.

(3)  $hg=g^{-1}h$, $gh\neq hg$, $g^{2}\neq h^{2}$,
$h_{1}g=g^{-1}h_{1}$, $gh_{1}\neq h_{1}g$ and
$gh_{1}\neq h_{1}g$. Lemma~\ref{ca0} $(ii)$, applied
to the elements $gx$ and $h$, gives us that either
$gxh=hgx=g^{-1}hx$ (and hence $g=g^{-1}$, a
contradiction), or $hgx=xg^{-1}h=xhg=hxg$ (and hence
$gx=xg$) or $(gx)^2=h^2$. So, we may assume that
$(gx)^2=h^2$. Similarly we that that
$(gx)^{2}=h_{1}^{2}$. Thus $h^{2}=h_{1}^{2}$. So
$hh_{1}\in N$ has order $2$ and thus is central in
$N$. However, this is impossible as $(h,h_{1})\neq
1$. This finishes the proof of the first part of the
statement.

Since $G_{[16,3]}$ contains noncentral elements of
order $2$, it thus follows at once that $G_{[16,3]}$
is not a subgroup of $G$. Since $N$ is not abelian
and elements of order $2$ are central, it hence
follows, from Lemma~\ref{ca0} $(i)$, that $\Char
(R)=4$ and that every nonabelian subgroup of $N$
generated by two elements is isomorphic with
$Q_{8}$. Hence, all subgroups of $N$ are normal in
$N$, i.e.  $N$ is a Hamiltonian $2$-group.
\end{proof}

\begin{trm}\label{n16,3Nab1}
Let $R$ be a commutative ring with $\Char(R) \neq 2$
and let $G$ be a nonabelian group of exponent 4 with
a nontrivial orientation homomorphism $\sigma$.
Assume that $N$ is  not abelian and that the
elements of order $2$ in $N$ are central in $N$.
Then $(RG)^-_{\varphi_{\sigma}}$ is commutative if
and only if $\Char(R)=4$ and  one of the following
conditions holds.
\begin{enumerate}
\item[$(i)$]   $G$ and $N$ are Hamiltonian 2-groups.
\item[$(ii)$] $G=\langle
a,b,c\ |\ a^4=c^4=1,\; b^2=a^2, ac=ca, ab=ba^{-1},
cb=bc^{-1} \rangle \times E$ and $N$ is equal to
either $\GEN{a,b}\times \GEN{c^2}\times E$ or
$\GEN{a,cb}\times \GEN{c^2}\times E$.
\end{enumerate}
\end{trm}

\begin{proof}  Assume $(RG)^-_{\varphi_{\sigma}}$
is commutative. Because of  Lemma~\ref{ord2cen} we
know  that the elements of order $2$ in $G$ are
central, $\Char(R)=4$, $G_{[16,3]}$ is not a
subgroup of $G$ and $N$ is a Hamiltonian $2$-group,
that is, $N=\GEN{a,b}\times E$, where
$\GEN{a,b}=Q_8$ and $E^2=1$.

Clearly, if $G$ contains an element $c$  of order
$2$ that does not belong to $N$, then $G=Q_{8}\times
E \times \langle c \rangle$. Hence $G$ also is a
Hamiltonian 2-group and we are in case (i) of the
theorem. So, to prove the necessity of the
conditions, we may suppose that the elements in
$G\setminus N$ have order $4$. As $E$ is central in
$G$, Remark~\ref{cenpro}  yields that we then have
that $G$ is the central product of $E$ and $\langle
a,b,c \rangle$, where $c\in G\setminus N$.  Hence
$G=\GEN{a,b,c}\times E_{1}$ for some subgroup
$E_{1}$ of $E$. Replacing, if necessary, $c$ by
either $ac$ or $bc$ we may assume that $c$ is not
central in $G$. So $(a,c)\neq 1$ or $(b,c)\neq 1$.

Next we show that $\GEN{a,b,c} =G_{[32,35]}$. By
Lemma~\ref{ca0} $(ii)$ we have that if $x\in
\GEN{a,b}$ with $\circ (x)=4$ then one of following
holds: $xc=cx$, $xc=c^{-1}x$ or $x^2=c^2$.

Assume there exists $x\in \GEN{a,b}$ with $\circ
(x)=4$ and $xc=cx$. Without loss of generality we
may assume that $x=a$. So $ac=ca$ and thus $bc\neq
cb$.  Note that then $a^{2}\neq c^{2}$ for otherwise
$ac$ is an element of order $2$ not contained in
$N$. It follows that $|\GEN{a,c}|=16$ and
$\GEN{a,c}=\GEN{a}\times \GEN{c}$. Clearly
$\GEN{a,c}\cap N=\GEN{a}\times \GEN{c^{2}}$ and,
since $\circ (b)=\circ (ab)=\circ (a^{-1}b) =4$, we
thus get that $b\not\in \GEN{a,c}$. Hence
$\GEN{a,b,c}\geq 32$. Note that, by the above,
$bc\neq cb$ and $b^{2}\neq c^{2}$ imply that
$bc=c^{-1}b$ then $\GEN{a,b,c}=\GEN{a,b,c^{2}}\cup
\GEN{a,b,c^{2}}c$ and $|\GEN{a,b,c^{2}}|=16$. So,
$|\GEN{a,b,c}|=32$ and it is easily seen that
$\GEN{a,b,c}=G_{[32,35]}$.

If $ac=c^{-1}a$ and $bc=c^{-1}b$ then $(ab)c=c(ab)$.
Because $\circ (ab)=4$ the previous yields that
$\GEN{a,b,c}=G_{[32,35]}$.

If $ac=c^{-1}a$ and $a^{2}=b^{2}=c^{2}$ then
$(bc)a=bac^{-1}=a^{-1}bc^{-1}=aa^{2}bc^{2}c=a(bc)$
and $b(bc)=c^{3}=c^{-1}b^{-1}b=(bc)^{-1}b$. Note
that $(bc)^{2}\neq 1$ as $bc\not\in N$. Hence $bc$
is not central. So, replacing $c$ by $bc$, we are
again in a situation that $c$ commutes with an
element of order $4$ in $\GEN{a,b}$. Hence,
$\GEN{a,b,c}=G_{[32,35]}$. The case $bc=c^{-1}b$ and
$a^{2}=b^{2}=c^{2}$ is dealt with similarly.

In order to finish the proof of the claim  we now
show that the following situation can not occur:
$a^{2}=b^{2}=c^{2}$, $bc\neq cb$, $ac\neq ca$,
$bc\neq c^{-1}b$ and $ac\neq c^{-1}a$.
Lemma~\ref{ca0} $(iii)$, applied to the elements
$ac$ and $bc$, yields that either
$(ac)(bc)=(bc)(ac)$, or $acbc = c^{-1}b^{-1}ac =
cc^{2} b^{2}b ac =  ca^{-1}bc = c^{-1}abc$ (and
hence $ac=c^{-1}a$, a contradiction) or
$acbc=bcc^{-1}a=ba=ab^{-1}$ (and hence
$cb=b^{-1}c^{-1}=bc$, a contradiction). Therefore
$(ac)(bc)=(bc)(ac)$ (i.e., $acb=bca$) and thus all
elements not in $N$ commute. In particular,
$(ac,acb)=1$ and thus $(ac,b)=1$. But then
$bca=acb=bac$ and thus $ca=ac$, again a
contradiction.

We already know that $G=\GEN{a,b,c}\times E_{1}$ and
$N=\GEN{a,b}\times E$ with $E_{1}$ a subgroup of the
elementary abelian $2$-group $E$.  Since
$\GEN{a,b,c}\cap N$ has index $2$ in $\GEN{a,b,c}$
and $c^{2}\not\in \GEN{a,b}$ it follows that
$E=\GEN{c^{2}}\times E_{1}$. Hence, the necessity of
the conditions follows from
Proposition~\ref{e32/35}.

The proof of the sufficiency follows from
Remark~\ref{obsexp} and Propositions~\ref{pca1} and
\ref{e32/35}.
\end{proof}

Now it is only left to classify the groups $G$ and
the kernels $N$ for which the
$\varphi_{\sigma}$-antisymmetric elements commute in
case
 $N$ contains a noncentral element of order
$2$. Then, by Lemma~\ref{fixedcentral}, we have that
$R_2=\{0\}$.
 In order to proceed with this case we first
prove the following lemma.

Assume $(RG)^-_{\varphi_{\sigma}}$ is commutative.
Recall from Lemma~\ref{ca2} $(i)$ that if $x$ and
$y$ are noncommuting elements of $N$ with $x^{2}\neq
1$ and $y^{2}=1$  then $\GEN{x,y}$ is either $D_{4}$
or $G_{[16,3]}$. In the next lemma we investigate
the structure of the group $\GEN{x,y,z}$ for $z\in
G\setminus N$.

\begin{lem}\label{add}
Assume $(RG)^-_{\varphi_{\sigma}}$ is commutative
and $G$ has exponent $4$.  Suppose $x,y\in N$ and
$z\in G\setminus N$ are such that $x^2\neq 1 = y^2$
and $z^2\neq 1$. If $xy\neq yx$ then $R_2=\{0\}$ and
$\GEN{x,y}$ is either $D_{4}$ or $G_{[16,3]}$.
Moreover
\begin{enumerate}
\item
If $\langle x,y\rangle=D_4$ then one of the
following conditions holds.
\begin{enumerate}
\item[(i)] $yz=zy$ and $x^2=z^2$;
\item[(ii)] $\langle y,z\rangle=D_4,$ and $xz=z^{-1}x=zx^{-1}$;
\item[(iii)] $\langle y,z\rangle=G_{[16,3]}$, $xz=zx$ and $x^2=(yz)^2$;
\item[(iv)] $\langle y,z\rangle=G_{[16,3]}$, $xz=z^{-1}x$ and $x^2=z^2$;
\item[(v)] $\langle y,z\rangle=G_{[16,3]}$, $xz=z^{-1}x$ and $x^2=(yz)^2$;
\item[(vi)] $\langle y,z\rangle=G_{[16,3]}$, $xz=z^{-1}x$ and $x^2=(yz)^2z^2$;
\item[(vii)] $\langle y,z\rangle=G_{[16,3]}$, $xyz=zxy$ and $x^2=z^2$.
\end{enumerate}
\item
If $\langle x,y\rangle=G_{[16,3]}$ then one of the
following conditions holds.
\begin{enumerate}
\item[(i)] $yz=zy$, $xz=zx$ and $z^2=x^2(xy)^2$;
\item[(ii)] $yz=zy$, $xz=z^{-1}x$ and $z^2=x^2(xy)^2$;
\item[(iii)] $\langle y,z\rangle=G_{[16,3]}$, $xz=zx$, $xyz=zyx$ and $z^2=(xy)^2$;
\item[(iv)] $\langle y,z\rangle=G_{[16,3]}$, $zxy=xyz=yzx$ and
$x^2=z^2$.
\end{enumerate}
\end{enumerate}
\end{lem}
\begin{proof} Suppose $xy\neq yx$. As mentioned before the
Lemma we already know that $R_2=\{0\}$ and
$\GEN{x,y}$ is either $D_{4}$ or $G_{[16,3]}$.
Since, by assumption, $G$ has exponent $4$,
Lemma~\ref{ca2} $(iii)$ yields  that $\langle
y,z\rangle$ is either abelian, $D_4$ or
$G_{[16,3]}$. Because of Lemma~\ref{ca0} $(ii)$, we
also have that
\begin{eqnarray} \label{mogel}
  xz=zx &  xz=z^{-1}x & \mbox{  or } x^2=z^2.
\end{eqnarray}

First, assume that $\GEN{x,y}=D_{4}$. If $\langle
y,z\rangle=D_4$ then, since $y^{2}=1$ and $z^{2}\neq
1$, we get that  $(yz)^2=1$. Therefore,
Lemma~\ref{ca2} $(ii)$, applied to $x$ and $yz$,
yields that $xyz=yzx$. Hence, $yx^{-1}z=yzx$ and
thus $x^{-1}z=zx$. Because $x^{2}\neq 1$, this
implies in particular that $xz\neq zx$. Hence,
(\ref{mogel}) yields that $xz=z^{-1}z$ or
$x^{2}=z^{2}$.  This with $x^{-1}z=zx$ implies that
$zx^{-1}=xz=z^{-1}x$. Therefore, $(ii)$ holds if
$\GEN{y,z}=D_{4}$. Assume now that $yz=zy$ (and thus
$(yz)^{2}\neq 1$). We claim that then $x^2=z^2$, and
thus (i) holds. Because of (\ref{mogel}), we may
assume that \begin{eqnarray} \label{mogel1}
 xz=zx &\mbox{ or } & xz=z^{-1}x.
\end{eqnarray}
Lemma~\ref{ca0} $(ii)$, applied to $x$ and $yz$,
yields that $xyz=yzx$, $xyz=z^{-1}yx$ or
$x^2=(yz)^2$. As $\GEN{x,y}=D_{4}$ this implies that
$x^{-1}z=zx$, $ x^{-1}z=z^{-1}x$ or $x^{2}=z^{2}$.
Because $x^{2}\neq 1$, we then obtain from
(\ref{mogel1}) that $x^{2}=z^{2}$.

Suppose now that $\langle y,z\rangle=G_{[16,3]}$.
Lemma~\ref{ca0} $(ii)$, applied  to $x$ and $yz$,
yields  that either $xyz=yzx$ (and hence
$x^{-1}z=zx$), $xyz=z^{-1}x^{-1}y$ or $x^2=(yz)^2$.
Because of (\ref{mogel}) we also know that $xz=zx$,
$xz=z^{-1}x$ or $x^{2}=z^{2}$. If $xz=zx$ then it
follows that $x^2=(yz)^2$ and therefore $(iii)$
holds. If $xz=z^{-1}x$ then we get that either
$x^2=z^2$, $(yz)^2z^2$ or $(yz)^2$, and therefore
$(iv)$, $(v)$  or $(vi)$ holds. If $x^2=z^2$ then
either $xz=z^{-1}x$ or $xyz=zxy$ and therefore
$(iv)$ holds or $(vii)$ holds.

Second, assume that $\GEN{x,y}=G_{[16,3]}$. First
assume that $\langle y,z\rangle=D_4$. Then
$(yz)^2=1$ and applying Lemma~\ref{ca2} $(ii)$, to
$x$ and $yz$ we have that $xyz=yzx$. Thus $xz\neq
zx$.  Since $(xy)^{2}\neq 1$, Lemma~\ref{ca2},
applied to $xy$ and $yz$, yields that $xyyz=yzxz$.
By Remark~\ref{remarkcentral} the element $y^{2}$ is
central. Hence it follows that $xzy=yzx$. This on
its turn implies that $xzy=z^{-1}yx$ and therefore
$xz\neq z^{-1}x$. From (\ref{mogel}) we thus get
that  $x^2=z^2$. Lemma~\ref{ca0} $(ii)$, applied to
$xy$ and $z$, also yields us that $zxy=xyz=xz^{-1}y$
(and thus $zx=xz^{-1}$) , $z^{-1}xy=xyz=xz^{-1}y$
(and thus $z^{-1}x=xz^{-1}$) or
$(xy)^{2}=z^{2}=x^{2}$ (and thus $yxy=x$). So each
time we obtain a contradiction. So, $\GEN{y,z}$ is
not $D_{4}$.

Assume now that $yz=zy$. Then $(yz)^2=z^2\neq 1$ and
applying Lemma~\ref{ca0} $(ii)$, to $x$ and $yz$, we
have that either $xyz=yzx$, $xyz=yz^{-1}x$ or
$x^2=z^2$. Thus, since we know that either $xz=zx$,
$xz=z^{-1}x$ or $x^2=z^2$, it follows that $z^2=x^2$
or $z^{2}=x^2(xy)^2$. If $z^2=x^2$ then, applying
Lemma~\ref{ca0} $(ii)$ to $yx$ and $z$ we get that
either $yxz=zyx=yzx$, $yxz=yz^{-1}x$ or
$(yx)^2=x^2$. Since $xy\neq yx$ we thus obtain that
$xz=zx$, or $xz=z^{-1}x$.  Applying Lemma~\ref{ca0}
$(ii)$ to $xy$ and $yz$, we get that either
$xz=yzxy$, $xz=yz^{-1}xy$ or $(xy)^2=x^2$. Again
because $yx\neq xy$, we obtain that $z^2=x^2(xy)^2$
and thus  $(xy)^2=1$, a contradiction. Therefore
$z^2=x^2(xy)^2$ and hence $(i)$ holds or $(ii)$
holds.

Finally assume that $\langle y,z\rangle=G_{[16,3]}$.
Lemma~\ref{ca0} $(ii)$, applied to $xy$ and $z$,
yields that either $xyz=zxy$, $xyz=z^{-1}xy$ or
$z^2=(xy)^2$. If $xyz=zxy$ or $xyz=z^{-1}xy$ then,
since $zy\neq yz$ and $yz\neq z^{-1}y$ and because
of (\ref{mogel}), we get that $x^{2}=z^{2}$. If
$z^2=(xy)^2$ then, since $xy\neq yx$, we get that
$z^2\neq x^2$. Hence, (\ref{mogel}) implies that
$xz= zx$ or $xz= z^{-1}x$. We claim that then
$xyz=zyx$. Suppose the contrary. Then,
Lemma~\ref{ca0} $(ii)$ applied  to $yx$ and $yz$,
gives us that  $yxyz=z^{-1}x$ or
$(yz)^2=(yx)^2=(xy)^2=z^2$. The former (together
with $xz=zx$ or $xz=z^{-1}x$) implies that $yxyz =
xz^{-1}$ or $yxyx=xz$. However this leads to a
contradiction because it results in
$z^2=x^2(xy)^2=x^2z^2$ and thus $x^{2}=1$. The
latter gives a contradiction as it implies $yz=zy$.

So we are left deal with two cases: (Case 1)
$xyz=zyx$ and $z^{2}=(xy)^{2}$, and (Case 2)
$x^{2}=z^{2}$, and $xyz=zxy$ or $xyz=z^{-1}xy$.

(Case 1): $xyz=zyx$ and $z^2=(xy)^2$. We show that
then $xz=zx$ and thus $(iii)$ holds. To prove this,
we apply  Lemma~\ref{ca0} $(ii)$ to $x$ and $yz$.
This yields that either $xyz=yzx$, $xyz=z^{-1}yx$ or
$x^2=(yz)^2$. If $xyz=zyx=yzx$ or $xyz=zyx=z^{-1}yx$
then either $yz=zy$ or $(zy)^2=1$, a contradiction.
Hence $x^2=(yz)^2$ and thus $xz^{-1}=x^2xyzy$. Since
$xyz=zyx$ this yields $xz^{-1}=x^{2}zyxy$. As
$x^{2}$ and $z^{2}$ are central, we obtain that
$xz^{-1} =zx(xy)^{2}= zx z^2=z^{-1}x$. Therefore
$xz=zx$, as claimed.

%
%

(Case 2): $x^2=z^2$, and $xyz=zxy$ or
$xyz=z^{-1}xy$. We will  prove that then
$zxy=xyz=yzx$ and thus $(iv)$ holds.

Clearly, $x^2\neq (yz)^2$. Hence, it follows from
Lemma~\ref{ca0} $(ii)$ that $xyz=yzx$ or
$xyz=z^{-1}yx$.

First assume that $xyz=yzx$. If $xyz=zxy$ then we
are done. If, on the other hand,  $xyz=z^{-1}xy$
then $yzx=xyz=z^{-1}xy$. Applying Lemma~\ref{ca0}
$(ii)$ to $yx$ and $yz$, we get that either
$yxyz=yzyx$, $yxyz=z^{-1}x$ or $(yx)^2=(yz)^2$.
Therefore we have that either $xyz=zyx$ and hence
$yz=zy$ a contradiction, $zx=z^{-1}x$ and hence
$z^2=1$, a contradiction or $xyx=zyz$ and hence
$yz=z^{-1}xyx=xyzx=x^2yz$,obtaining that $x^2=1$, a
contradiction.

Second assume that assume that $xyz=z^{-1}yx$. If
$xyz=z^{-1}xy$ we have that $xy=yx$, a
contradiction. Therefore to end the proof of the
lemma we have to deal with the case that $x^2=z^2$
and $zxy=xyz=z^{-1}yx$. Applying Lemma~\ref{ca0}
$(ii)$ to $yx$ and $yz$, we get that either
$yxyz=yzyx$, $yxyz=z^{-1}x$ or $(yx)^2=(yz)^2$. In
the first case we have that $xyz=zyx$ and since
$xyz=zxy$ we have that $(x,y)=1$, a contradiction.
In the second case we have that
$z^{-1}x=yxyz=yz^{-1}yx$ and hence $(z,y)=1$, again
a contradiction. Finally if $(yx)^2=(yz)^2$ then
$xyx=zyz$ and hence $z^{-1}xyx=yz$. Since $zxy=xyz$
and $z^2=x^2$ is central it follows that
$yz=z^{-1}xyx=zxyx^{-1}=xyzx^{-1}$. Therefore we get
that $xyz=yzx$ as desired.
\end{proof}

\begin{trm}\label{FinalTheorem}
Let $R$ be a commutative ring with $\Char(R) \neq 2$
and let $G$ be a nonabelian group of exponent $4$
with a nontrivial orientation homomorphism $\sigma$.
Assume that $N$ is  not abelian and that there
exists a noncentral element of order $2$ in $N$.
Then $(RG)^-_{\varphi_{\sigma}}$ is commutative if
and only if $R_2=\{0\}$ and one the following
conditions holds
\begin{enumerate}
\item[$(i)$]   $G\cong\langle a,b, c \mid
a^2=b^2=c^2=1,\; abc=bca=cab \rangle\times E$ and
$N\cong D_4\times E$, where $E^2=1$;
\item[$(ii)$] $G=\langle a,b,c,d\, |\, a^4=b^2=c^2=d^2=1, ab=ba, ac=ca, ad=dab, bc=cb, bd=db, cd=da^2c
\rangle\times E$ and $N=\langle b\rangle\times
\langle c,d\rangle\times E$, where $E^2=1$;
\item[$(iii)$] $G=\langle a,b,c\, |\, a^4=b^4=c^2=1, ab=ba, ac=ca^{-1}, bc=ca^2b^{-1}\rangle\times E$ and $N=\langle
a,c\rangle\times \langle b^2\rangle\times E$, where
$E^2=1$.
 \item[$(iv)$] $G=\langle g,a,b\, |\, g^4=a^4=b^2=1,\;
ga=ag,\; gb=bg,\; ab=g^2ba\rangle\times E$ and
 $N=\langle a,b\rangle \times E$
 or $N=\GEN{ ga,b} \times E$, where $E^2=1$.
\end{enumerate}
\end{trm}
\begin{proof} The sufficiency of the conditions
follows from Remark~\ref{obsexp},
Proposition~\ref{pca4}, Proposition~\ref{e32/38},
Proposition~\ref{e32/39} and
Proposition~\ref{e32/16}.

To prove the necessity,  assume
$(RG)^-_{\varphi_{\sigma}}$ is commutative. Let
$N=\ker (\sigma )$. So $N$ has index $2$ in $G$.
Since by assumption $N$ contains a noncentral
element of order $2$, Lemma~\ref{fixedcentral}
 yields that $R_2=\{0\}$. Hence by Lemma~\ref{ca1} $(i)$ we get that $Q_8\not\subseteq N$.

 Since $N$ is not
abelian, it then follows from  \cite{BP} (see
Theorem~3.3 and the introduction of Section 4) that
either
\begin{enumerate}
\item[] (Case 1)  The elements of order $4$ in
$N$ generate an abelian subgroup.
\item[]
(Case 2) $N$ contains an elementary abelian
$2$-group of index 2.
\end{enumerate}

(Case 1)  Assume  $A=\GEN{x\in N \mid x^2\neq 1}$ is
an abelian subgroup. It follows from \cite{JR2} that
$A$ is a subgroup of index $2$ in $N$.
Thus, write $G=N\cup Nh$ and $N=A\cup Ay$ for some
$h\in G$ and $y\in N$ with $y^{2}=1$. Since
$G_{[16,3]}$ is a nonabelian group generated by
elements of order $4$, it follows that it is not
contained in $N$. Furthermore, as $N$ is not
abelian, Lemma~\ref{ca2} $(i)$ yields that for $a\in
A$, either $(a,y)=1$, or $\GEN{a,y}\cong D_{4}$ and
$(ay)^{2}=1$. Hence, we can choose $x\in A$, with
$x^{2}\neq 1$ and so that $\GEN{x,y}=D_{4}$. We also
note that we may assume that $h^{2}\neq 1$. Indeed,
suppose $h^{2}=1$. Then, again Lemma~\ref{ca2}
$(ii)$, $(h,x)=1$ and thus
$(hx)^{2}=h^{2}x^{2}=x^{2}\neq 1$. So, replacing $h$
by $hx$ if needed, we indeed may assume that
$h^{2}\neq 1$. Hence, by Lemma~\ref{ca2} (iii),
$\GEN{y,h}=D_{4}$ and thus $yh=h^{-1}y$.

Let $g\in N$  with $g^2\neq 1$ (thus $g\in A$).
Since elements of $N\setminus A$ have order $2$, we
get that $(y,g)\neq 1$. Furthermore, by the above,
$\GEN{y,g}=D_{4}$ and thus $yg=g^{-1}y$. By
Lemma~\ref{ca0} ($ii$) we have that $hg=gh$,
$gh=h^{-1}g$ or $h^2=g^2$. We claim that
 \begin{eqnarray} g^2=h^2 &  g^2=(yh)^2 &  \mbox{ or }
 g^2=h^2(yh)^2.\label{claimmogel}
 \end{eqnarray}
First, assume  that $hg=gh$. Then $(yh,g)\neq 1$
and, by Lemma~\ref{ca2} ($ii$), we have that
$(yh)^2\neq 1$. Therefore, Lemma~\ref{ca0} ($ii$),
yields that either $gyh=h^{-1}yg$  or $(yh)^2=g^2$
as desired. The latter is as desired. In the former
case, $ghyh=hgyh=yg$ and thus $g^{-1}yhyh=g$. Hence
$(yh)^{2}=g^{2}$, again as desired in the claim.
Second, assume  that $gh=h^{-1}g$. Let
$K=\GEN{yh,g}$. If $K$ is abelian then
$yhg=gyh=g^{2}ygh=g^{2}yh^{-1}g=g^{2}yh^{2}hg=g^2
h^2 yhg$ and thus  $h^2=g^2$, again  as desired. If
$K$ is not abelian then, by Lemma~\ref{ca2} ($ii$),
$(yh)^2\neq 1$.  Lemma~\ref{ca0} therefore yields
that $g^{2}=(yh)^{2}$ or
$gyh=h^{-1}yg=h^{-1}g^{-1}y= g^{-1}hy$ and thus
$g^{2}=hyh^{-1}y^{-1}=(hy)^{2}h^{2}=h^{2}(yh)^{2}$,
again as desired. This proves the claim
(\ref{claimmogel}).

We now prove the following five statements.
\begin{enumerate}
\item[] (1.a): $\z
(N)=\{ a\in A \mid a^{2}=1\}$
and $A/\z (N)$ is an elementary abelian $2$-group.
\item[] (1.b): $A/\z (N)$ is cyclic, or
equivalently, $A=\GEN{x} \times E_{1}$ for some
elementary abelian $2$-group $E_{1}$.
\item[] (1.c): $\z (N)\subseteq \z (G)$.
\item[] (1.d): $G=\GEN{x,y,h}\times E$ for
some elementary abelian $2$-subgroup $E$ of $G$.
\item[] (1.e): $\GEN{x,y,h}$ is isomorphic with
either $G_{[32,30]}$, $G_{[16,13]}$ or
$G_{[32,31]}$.
\end{enumerate}
It then follows from Remark~\ref{obsexp},
Proposition~\ref{pca4}, Proposition~\ref{e32/38} and
Proposition~\ref{e32/39} that either condition (i),
(ii) or (iii) of the statement of the result is
satisfied. This then finishes the proof of (Case 1).

\vspace{12pt} (1.a) First we show that if $1\neq
z\in \z (N)$ then $z$ has order $2$. Indeed, for
suppose $z$ has order $4$. Then $z\in A$. Since
$y^{2}=1$ it follows that $yz$ has order $4$. Hence
also $yz\in A$, a contradiction.

Since $A$ is abelian, $N=A\cup Ay$ and $(y,x)\neq
1$, it follows that $\z (N)\subseteq A$. Hence, $\z
(N) \subseteq \{ a\in A \mid a^{2}=1\}$.
Conversely, let $a\in A$ with $a^{2}=1$. If
$(y,a)\neq 1$ then,  by Lemma~\ref{ca3},
$\GEN{y,a}=D_4$. Hence, $(ya)^2\neq 1$ and $ya\in
N\setminus A$, a contradiction. So, $(y,a)=1$ and
thus $a\in \z (N)$. So we have shown that $\z (N)
=\{ a\in \mid a^{2}=1\}$.

Because of Remark~\ref{remarkcentral}, we also know
that squares of elements of $G$ are central. In
particular, $A/\z (N)$ is an elementary abelian
$2$-group.

\vspace{12pt} (1.b) Because of part (1.a), in order
to prove this property, it is sufficient to show
that there does not exist an element $g\in A$ of
order $4$ so that $\GEN{g,x} = \GEN{g} \times
\GEN{x}$. Assume the contrary. By
(\ref{claimmogel}), we know that $$x^{2}=h^{2},\;\;
x^{2}=(yh)^{2} \mbox{ or } x^{2}=h^{2}(yh)^{2}.$$ We
will show that each of these cases leads to a
contradiction. Note that, also by
(\ref{claimmogel}), $g^{2}=h^{2}$, $g^{2}=(yh)^{2}$
or $g^{2}=h^{2}(yh)^{2}$.

Assume that $x^2=h^2$. Since, by assumption
$g^{2}\neq x^{2}$, Lemma~\ref{ca0} ($ii$) yields
that we $gh=hg$ or $gh=h^{-1}g$.

Suppose  $gh=hg$. Since $(y,g)\neq 1$ (see above),
we have that $(yh,g)\neq 1$ and therefore, by
Lemma~\ref{ca2} ($ii$), $(yh)^2\neq 1$. Applying
Lemma~\ref{ca0} ($ii$) to $yh$ and $g$, we deduce
that  $gyh=h^{-1}yg$ (and hence $g^2=(yh)^2$) or
$g^2=(yh)^2$. So, $g^2=(yh)^2$ and hence (as squares
are central) $(y,h)=h^2(yh)^2=h^2g^2$. Now, applying
Lemma~\ref{ca0} ($ii$) to the elements $yh$ and $x$,
we get that  $yhx=xyh$, $xyh=h^{-1}yx$ or
$x^2=(yh)^2$. We now show that each of these three
cases leads to a contradiction.
 If $yhx=xyh=yx^{-1}h$ then $hx=x^{-1}h$ and hence
$(x,h)=x^2=h^2$. On the other hand,
$(yxh)^2=yxhyxh=(yh)^2x^2=x^2g^2\neq 1$, because of
the assumption. Since also $(g,yhx)\neq 1$,
Lemma~\ref{ca0} ($ii$) therefore implies that
$gyxh=h^{-1}x^{-1}yg=hxyg=g^{-1}hxy$. Hence,
$g^2=(yxh)^2=x^2g^2$ and thus  $x^2=1$, a
contradiction.
 If $xyh=h^{-1}yx=h^{-1}x^{-1}y=hxy$, then
since $xyh=x(y,h)hy=xy^{2}h^{2}(yh)^{2}hy=
h^2g^2xhy$, we get that  $h^2g^2xh=hx$. Therefore,
$h^{2}g^{2}=(h,x)$. We also know that
$(h,x)=h^{2}x^{2}(hx)^{2}=(hx)^{2}$. Thus
$h^{2}g^{2}=(hx)^{2}$. Then, consider the group
$\GEN{gx,h}$. By Lemma~\ref{ca0} ($ii$), we get that
either $gxh=hgx=ghx$ (and hence $xh=hx$, a
contradiction), or $gxh=h^{-1}gx=gh^{-1}x$ (and
hence $xh=h^{-1}x$; so that $h^{2}g^{2}=
(xh)^{2}=x^{2}$ and thus, because $h^{2}=x^{2}$, we
get that $g^{2}=1$, a contradiction), or
$(gx)^{2}=h^{2}$ (and hence$g^{2}x^{2}=h^{2}$ and
therefore $g^{2}=1$, a contradiction). So
$xyh=h^{-1}yx$ is excluded.
 If $x^2=(yh)^2$, then,  since
$(yh)^2=g^2$, we  get $x^{2}=g^{2}$, again  a
contradiction. This shows that if $x^{2}=h^{2}$ then
$gh\neq hg$.

 Therefore $x^{2}=h^{2}$ implies
that $gh=h^{-1}g$. Notice that $(h,x)\neq 1$ because
otherwise $(hx)^2=1$ and hence, by Lemma~\ref{ca2}
($ii$), $ghx=hxg=hgx$. Then $hg=gh=h^{-1}g$ and
hence $h=h^{-1}$, a contradiction. So, $(gh,x)\neq
1$ and, since $(gh)^{2}=g^{2}\neq x^{2}$, applying
Lemma~\ref{ca0} ($ii$) to $gh$ and $x$, one deduces
that $xgh=h^{-1}g^{-1}x=h^{-1}xg^{-1}$. Since
$xgh=xh^{-1}g$ and $h^{2}=x^{2}$, we get that
$g^2=hx^{-1}h^{-1}x=(hx)^2$ and thus
$(x,h)=(hx)^2=g^2$. It follows that
$(gx,h)=g^{2}h^{2}\neq 1$, because otherwise
$g^{2}=h^{2}=x^{2}$, a contradiction. Therefore, by
Lemma~\ref{ca0} ($ii$), $gxh=h^{-1}gx$ or
$(gx)^{2}=h^{2}$. The former implies that
$g^{2}h^{2}hgx=h^{-1}gx$ and thus $g^{2}=1$,  a
contradiction. The latter yields that
$g^{2}x^{2}=(gx)^{2}=h^{2}=x^{2}$ and thus
$g^{2}=1$, again a contradiction.

So we have shown that indeed $x^{2}\neq h^{2}$.
Since $x$ and $g$ play a symmetric role and
$x^{2}\neq g^{2}$, we also get that $h^{2}\neq
g^{2}$ and it only remains to deal with the case
that $x^{2}=(yh)^{2}$ and $g^{2}=h^{2}(yh)^{2}$ (and
thus $(yh)^{2}=x^{2}=g^{2}h^{2}$).
%
Again by Lemma~\ref{ca0} ($ii$), we have that
$gh=hg$ or $gh=h^{-1}g$. Assume first that
$(g,h)=1$. Notice that then $(yh,g)\neq 1$ and
hence, by Lemma~\ref{ca2} ($ii$), $(yh)^2\neq 1$.
Then, by Lemma~\ref{ca0} ($ii$) applied to $yh$ and
$g$, we get that  $g^2=(yh)^2$ or
$gyh=h^{-1}yg=g^{-1}h^{-1}y$. In both cases this
implies $g^2=(yh)^2$. Then $h^2(yh)^2=g^2=(yh)^2$
and therefore $h^2=1$, a contradiction. Thus,
$gh=h^{-1}g$. If $(hg,x)\neq 1$ then, by
Lemma~\ref{ca0} ($ii$), we have that
$g^2=(hg)^2=x^2$, a contradiction, or
$xhg=g^{-1}h^{-1}x=g^2h^2ghx=g^2h^2h^{-1}gx=g^2hxg$.
So, $g^{2}=(h,x)$. Hence,
$g^{2}=(h,x)=h^{2}x^{2}(hx)^{2}=g^{2}(hx)^{2}$ and
thus $(hx)^{2}=1$.
Then, by Lemma~\ref{ca2} ($ii$), we have that
$hxg=ghx=h^{-1}gx=h^{-1}xg$ and thus $h^{2}=1$, a
contradiction. Hence,
$(hg,x)=1$. Then $(h,x)=1$. Notice that
$(hx)^{2}\neq 1$. Indeed, for otherwise, by
Lemma~\ref{ca2} ($ii$), $(hx,x)=1$ and thus
$(h,x)=1$, a contradiction. So,
$1=(hx)^{2}=h^{2}x^{2}$ and thus $h^{2}=x^{2}$, a
contradiction. Also $(hg)^{2}=g^{2}\neq 1$.
Applying Lemma~\ref{ca0} ($iii$) to $hx$ and $hg$,
we get that  $hxhg=hghx$ (and hence $hg=gh$, a
contradiction), or $hxhg=g^{-1}h^{-1}hx$ (and thus
$h^2=g^2$, a contradiction),  or
$hxhg=hgx^{-1}h^{-1}$ (and hence $x=x^{-1}$, a
contradiction) or $1=hxhghxhg=x^2h^2h^2g^2$ (and
hence $x^2=g^2$ again a contradiction).

So this finishes the proof of (1.b).

\vspace{12pt} (1.c)  We prove that $\z (N) \subseteq
\z (G)$.
%
So, let $e\in \z(N)$. By means of contradiction
assume that $(h,e)\neq 1$. If $(he)^2=1$ then, by
Lemma~\ref{ca2} ($ii$), we have that $hex=xhe$. As
$ex=xe$, we thus get that $hx=xh$. Also, by
Lemma~\ref{ca2} ($ii$), we get that $xehe=hexe$.
Thus, $(xe)(he) =hxe^2=hx=xh=(xe)(eh)$ and hence
$he=eh$, a contradiction. So, we also may assume
that $(he)^2\neq 1$. By Lemma~\ref{ca0} (ii), we
have that $hx=xh$, $xh=h^{-1}x$ or $x^2=h^2$. We now
prove that each case leads to a contradiction.

First, assume  $hx=xh$. Then, $(h,xe)\neq 1$ and
hence, by Lemma~\ref{ca0} (ii), we have that either
$xeh=h^{-1}xe=xh^{-1}e$  or $h^{2}=x^{2}$. The
former yields $(he)^2=1$, a contradiction. The
latter implies $(hx)^2=1$ and thus, by
Lemma~\ref{ca2} ($ii$), we get that
$hxxe=xehx=ehx^2$. So, $he=eh$, a contradiction.

Second, assume  $xh=h^{-1}x$. Applying
Lemma~\ref{ca0} ($ii$) to $h$ and $xe$, we get that
$hxe=xeh=exh=eh^{-1}x$, or $xeh=h^{-1}xe=xhe$ or
$h^{2}=x^{2}$. The former leads to $(he)^2=1$, a
contradiction. The second implies $he=eh$, a
contradiction. So, $h^2=x^2$. Since $(x,h)\neq 1$,
we obtain that $(he,x)\neq 1$. Then, applying
Lemma~\ref{ca0} ($ii$) to $he$ and $x$, we get that
$xhe=eh^{-1}x=exh=xeh$ and hence $he=eh$, a
contradiction, or $(he)^2=x^2=h^2$ and hence
$eh=he$, a contradiction.

Third, assume  $x^2=h^2$, $(h,x)\neq 1$ and $xh\neq
h^{-1}x$. If $hexe=xehe$ then
$hxe=xeh=x^{3}eh^{3}=x^{-1}eh^{-1}=ex^{-1}h^{-1}$
and hence $(hxe)^2=1$. Then, by Lemma~\ref{ca2}
($ii$), we have that $(hxe,x)=1$ and hence
$(h,x)=1$, a contradiction. So, $(he,xe)\neq 1$.
Therefore, applying Lemma~\ref{ca0} ($ii$) to $he$
and $xe$, we get that $xehe=eh^{-1}xe$ and hence
$xh=h^{-1}x$, a contradiction, or $(he)^2=x^2=h^2$
and hence $eh=he$, again a contradiction. This
finishes the proof of  (1.c).

\vspace{12pt} (1.d)  This follows at once from
(1.b), (1.c) and Remark~\ref{cenpro}.

\vspace{12pt} (1.e) We determine the group
$\GEN{x,y,h}$. Recall that either $x^{2}=h^{2}$,
$x^{2}=(yh)^{2}$ or $x^{2}=(yh)^{2}h^{2}$ (see
(\ref{claimmogel})). Also, remember that
$\GEN{x,y}=D_{4}$, $y^{2}=1$, $\circ (h)=4$, and
thus, because of Lemma~\ref{ca2} ($iii$), either
(1.e.i) $yh=hy$, (1.e.ii) $yh=h^{-1}y$ or (1.e.iii)
$\GEN{y,h}=G_{[16,3]}$. Note that $|\GEN{x,y,h}|\leq
32$. We will deal with each of the three cases
separately.

(1.e.i) Suppose $yh=hy$. Then $(yh)^2=h^2$. Since
$x^2=h^2$, $x^2=(yh)^2$, or $x^2=(yh)^2h^2$, we thus
get that $x^2=h^2$. Hence,
$(h,x)=h^{2}x^{2}(hx)^{2}=(hx)^{2}$ and
$(h,x)=(x,h)=x^{2}h^{2}(xh)^{2}$. So,
$(xh)^2=(xh)^{-2}=(h^{-1}x^{-1})^2=(hx)^2$.

If $(xh)^2=1$ then $\GEN{x,y,h}=\GEN{xy,y,xh}$.
Since $\circ (xy)=\circ (y)=\circ (xh)=o(c)$ and
$(xy)y(xh)= y(xh)(xy)=(xh)(xy)y =x^{2}h$, we obtain
that $\GEN{x,y,h}=G_{[16,13]}$, as desired.

If $(xh)^2\neq 1$ then  let $a=h$, $b=(xh)^2$, $c=y$
and $d=xy$. Clearly, $\circ (a)=4$, $\circ (b)=
\circ (c)=\circ (d)=2$,
$dab=xyh(xh)^2=yx^2(xh)^3=yx^2h^{-1}x^{-1}=yhx^{-1}=hxy=ad$,
$da^2c=xyh^2y=y^2xh^2=y^2x^{-1}=yxy=cd$, $ab=ba$,
$ac=ca$, $bc=cb$ and $bd=db$. It follows that
$\GEN{x,y,h}=G_{[32,30]}$, again as desired.

(1.e.ii) Suppose  $yh=h^{-1}y$. Then, $(yh)^{2}=1$
and thus,  by Lemma~\ref{ca2} $(ii)$, we have that
$xyh=yhx$. Hence, $xh=hx^{-1}$. Since, by
Lemma~\ref{ca1} ($ii$), $xh=hx$, $xh=h^{-1}x$ or
$x^2=h^2$, we have that $x^2=1$, $x^2=h^2$ or
$xh=h^{-1}x$, respectively. Therefore if
$yh=h^{-1}y$ then $xh=h^{-1}x$ and $x^2=h^2$. Then
$\GEN{x,y,h}=\GEN{xy,y,yh}$ with $\circ (xy) = \circ
(y) = \circ(yh) = 2 $ and
$(xy)y(yh)=y(yh)(xy)=(yh)(xy)y=xyh$. Thus
$\GEN{x,y,h}= G_{[16,13]}$, as desired.

(1.e.iii) Suppose $\GEN{y,h}=G_{[16,3]}=H\cup Hh$,
where $H=\GEN{y}\times\GEN{h^2}\times\GEN{(yh)^2}$
is an elementary abelian $2$-group of order $8$.
Recall from Lemma~\ref{ca0} $(ii)$  that $xh=hx$,
$xh=h^{-1}x$ or $x^2=h^2$. We deal with each of
these cases separately.

If $x^2=h^2$ then $|\GEN{x,y,h}|=32$. Since $h^2\neq
(yh)^2$,  Lemma~\ref{ca0} $(ii)$ yields that
$x(yh)=(yh)x$ or $xyh=(yh)^{-1}x$, and thus
$xh=h^{-1}x$ or $xyh=hxy$. If $xyh=hxy$ then
$\GEN{x,y,h}=G_{[32,30]}$. For this it is enough to
note that $\GEN{x,y,h}=\GEN{a,b,c}$, with $a=h$,
$b=(xh)^2$, $c=xy$ and $d=y$, and $\circ (a)=4$,
$\circ (b)= \circ (c) =\circ (d) =2$,
$dab=yh(xh)^2=yx^{-1}(xh)^{-1} = xyh^{-1}x^{-1} =
h^{-1}xyx^{-1} = hy = ad$,
$da^2c=yh^2xy=yxyh^2=x=cd$, $ab=ba$, $ac=ca$,
$bc=cb$ and $bd=db$. If $xh=h^{-1}x$ then
$\GEN{x,y,h}=G_{[32,31]}$. To see this it is enough
to note that $\GEN{x,y,h}=\GEN{a,b,c}$, with $a=x$,
$b=yh$, $c=y$, $\circ (b) = 4$,
$ab=xyh=yh^{-1}x^{-1}=ba$, $ac=xy=ca^{-1}$ and
$ca^2b^{-1}=yx^2h^{-1}y=yhy=bc$.

Suppose now that $xh=h^{-1}x$ and $x^{2}\neq h^{2}$.
From (\ref{claimmogel}) we know that we have to
consider three cases: $x^{2}=h^{2}$,
$x^{2}=(yh)^{2}$ or $x^{2}=(yh)^{2}h^{2}$. The
former of course is excluded. If $x^2=(yh)^2$ then
$\GEN{x,y,h}=\GEN{a,b,c,d}$, with $a=xh$, $b=h^2$,
$c=xy$, $d=y$, and $\circ (a) = 4$, $\circ (b) =
\circ (c) = \circ (d) = 2$,
$ad=xhy=h^{-1}yx^{-1}=yh(yh)^2x^{-1}=yhx=yxhh^2=dab$,
$da^2c=y(xh)^2xy=yx^{-1}y=x=cd$, $ab=ba$,
$ac=xhxy=x^2(yh)^{-1}=x^2(yh)^2yh=yh=xyxh=ca$,
$bc=cb$ and $bd=h^2y=yh^2=db$; so
$\GEN{x,y,h}=G_{[32,30]}$. If $x^2=(yh)^2h^2$ then
take $a=xh$, $b=h^2$, $c=y$ and $d=xy$. Because
$\circ (a) = 4$, $\circ (b)= \circ (c) =\circ (d) =
2$,
$ad=xhxy=x^2(yh)^{-1}=x^2(yh)^2yh=h^2yh=xyxhh^2=dab$,
$da^2c=xyx^2y=yxy=cd$, $ab=ba$,
$ac=xhy=(yh)^{-1}x^{-1}=yh(yh)^2h^2x=yhx^2x=yxh=ca$,
$bc=cb$ and $bd=db$, it follows that
$\GEN{x,y,h}=G_{[32,30]}$.

Assume that $xh=hx$. If $x^2\neq(yh)^2$ then, by
Lemma~\ref{ca0} $(ii)$, we have that $x(yh)=(yh)x$
or $x(yh)=(yh)^{-1}x$. Since $xyh=yhx^{-1}$ we get
that $x^2=1$ or $x^2=(yh)^2$, a contradiction. Thus,
$x^{2}=(yh)^{2}$. Let  $a=x$, $b=h$, $c=y$. Clearly,
$ab=ba$, $ac=ca^{-1}$ and
$ca^2b^{-1}=yx^2h^{-1}=y(yh)^2h^{-1}=hy=bc$. Hence,
since $|\GEN{x,y,h}|=32$, we obtain that
$\GEN{x,y,h}=G_{[32,31]}$.


This finishes the proof of (1.e) and hence also the
proof of (Case 1).

\vspace{12pt} (Case 2) Assume $N$ contains an
elementary abelian $2$-subgroup $B$ of index $2$ and
that the elements of order $4$ in $N$ do not
generate an abelian subgroup. We claim that if $c\in
N$ with $c^2\neq 1$ and  $a\in\setminus B$ with
$a^2=1$
 then $\GEN{a,c}$ is either abelian or $D_4$. Indeed,  assume that $(a,c)\neq 1$. Then
 by Lemma~\ref{ca2} then either
 $\GEN{a,c}=D_4$ or  $\GEN{a,c}=G_{[16,3]}$.  Since $c=ab$ for some $b\in B$ it follows that
 $\GEN{a,c}=\GEN{a,b}$, a contradiction, because  $G_{[16,3]}$ can not
 be generated by two elements of order 2. This proves the claim.

Next we claim that $a^2\neq 1$ for all $a\in
N\setminus B$. Suppose the contrary, then by the
previous claim we have that for all $c\in N$ with
$c^2\neq 1$, $\GEN{a,c}$ is either abelian or $D_4$.
By the assumptions there exist $b_1a, b_2a\in N$
both of order $4$ and $b_1,b_2\in B$ so that
$(b1_a,b_2a)\neq 1$. Since $(b_i,a)\neq 1$ ($i=1,2$)
it follows that $\GEN{a,b_i}=D_4=\GEN{b_ia,b_j}$
($j,i=1,2$). Hence $b_1a b_2
a=(ab_2)^{-1}b_1a=b_2ab_1a$, a contradiction. This
finishes the proof of the claim.

 As a consequence of the previous claim we have that $D_4$ can not be a subgroup
of $N$, because otherwise we can always find an
element of order $2$ in $N\setminus B$.

Since $Q_{8}$ is not contained in $G$,
Lemma~\ref{ca1}  yields that $N$ contains
$G_{[16,3]}$. So, let $c\in N$ with $c^2\neq 1$ and
$b\in B$ such that $\GEN{b,c}=G_{[16,3]}$.

 Assume there exists $g\in
G\setminus N$ with $g^{2}=1$. Then by
Lemma~\ref{ca2} $(ii)$ it follows that $(g,c)=1$ and
$1=(g,cb)$. Therefore $1=(g,b)$ and hence
$(gb)^2=1$. Again by  Lemma~\ref{ca2} $(ii)$ we have
that $1=(gb,c)=(c,b)$, a contradiction. So we have
shown that $g^2\neq 1$ for all  $g\in G\setminus N$.

Further, also choose $a\in N\setminus B$. Then,
because of the claims above and Lemma~\ref{ca2}
$(i)$, for all $b\in B$ with $(a,b)\neq 1$ it
follows that $$\GEN{a,b}=G_{[16,3]}.$$

Now we are going to show that $N=\GEN{a,b}\times E$
for  some elementary abelian $2$-subgroup $E$ of $\z
(N)$. Let $g\in G\setminus N$, so $g^2\neq 1$. First
we deal with the case $g^2=a^2$. If $b,b_1\in B$
such that $(a,b)\neq 1\neq (a,b_1)$, then by
Lemma~\ref{add} (2) it follows that
$(a,b)=(g,a)=(a,b_1)$. Therefore $(a,bb_1)=1$ and it
follows that $B=\GEN {b} \times \GEN{a^{2}}\times
E$, for some elementary abelian $2$-subgroup $E$ of
$\z (N)$. Thus $N=\GEN{a,b}\times E$ as desired.
Second we deal with the case $g^2\neq a^2$. If
$(g,a)\neq 1$, then by Lemma~\ref{add} (2) it
follows that $(a,b)=g^2=(a,b_1)$. Therefore
$(a,bb_1)=1$ and again as above we have that
$N=\GEN{a,b}\times E$ as desired. Finally if
$(g,a)=1$ then again by Lemma~\ref{add} (2) the
commutators $(a,b)$ and $(a,b_1)$ are either $g^2$
or $(ga)^2$. In case $(a,b)=(a,b_1)$ arguing as
before we obtain the desired conclusion. So assume
that $(a,b)\neq (a,b_1)$. But then, again by
Lemma~\ref{add} (2) we have that
$(a,bb_1)=(a,b)(a,b_1)=g^2(ga)^2$ is either $g^2$ or
$(ga)^2$. This is in contradiction
 with the fact that  the elements in $G\setminus N$ are all of order 4.

We now show that there exists $g\in G\setminus N$
such that
$$\GEN{g,a,b}=\GEN{g,a,b\,|\;g^4=a^4=b^2=1,
\;ga=ag,\;gb=bg,\;ab=g^2ba}=G_{[32,24]}.$$ For this
note that by Lemma~\ref{ca0} $(ii)$, for any $g\in
G\setminus N$ we have  either $(1)$ $ag=ga$ or $(2)$
$ag=g^{-1}a$ or $(3)$ $a^2=g^2$. First assume
$(1)$, that is $ag=ga$. Then it is easy to verify
that case (2) (i) or  case (2) (iii) of
Lemma~\ref{add} must holds. In the first case  it is
readily verified that $\GEN{g,a,b}=G_{[32,24]}$ and
in the second case, replacing $g$ by $ga$, one also
obtains that
$\GEN{g,a,b}=\GEN{ga,a,b}=G_{[32,24]}$. 
Assume now $(2)$, that is $ag=g^{-1}a$. Then (2)
$(ii)$ of Lemma~\ref{add} holds. Thus $gb\in
\z(\GEN{g,a,b})$ and $(gb)^2=(a,b)$. Therefore
$\GEN{g,a,b}=\GEN{gb,a,b}=G_{[32,24]}$. Third assume
$(3)$, that is,  $a^2=g^2$. Then (2) (iv) of
Lemma~\ref{add} holds. Thus $gab\in \z(\GEN{g,a,b})$
and $(gab)^2=(a,b)$. Therefore
$\GEN{g,a,b}=\GEN{gab,a,b}=G_{[32,24]}$ as desired.

Now we are going to prove that $E\subseteq \z(G)$
and therefore $G=G_{[32,24]}\times E$ finishing the
proof of the theorem. Let $e\in E$. We need to show
that $(g,e)=1$. From the above we know that there
exists $g\in G\setminus N$ such that $\GEN{g,a,b}=
 \GEN{g,a,b\,|\;g^4=a^4=b^2=1, \;ga=ag,\;gb=bg,\;ab=g^2ba}$.
  Applying
Lemma~\ref{ca0} $(ii)$ to $ae$ and $g$ we have that
either $ae g=gae$ and thus $(g,e)=1$ as desired; or
$aeg=g^{-1}ae$ then $eg=g^{-1}e$ and thus
$(ge)^2=1$, a contradiction because $ge\not \in N$;
or $g^2=(ae)^2$ and thus $g^2=a^2$, a contradiction.
\end{proof}

%
%
%
%

\begin{tabular}{lll}
O. Broche Cristo &Eric Jespers &  Manuel Ruiz
\\ Dep. de Ci\^{e}ncias Exatas &Dept. Mathematics  & Dep.
M\'{e}todos Cuantitativos e Inform\'{a}ticos \\
Univ. Federal de Lavras &Vrije Universiteit Brussel
&Univ. Polit\'{e}cnica de Cartagena
\\ Caixa Postal 3037& Pleinlaan 2 &  Paseo Alfonso
XIII, 50 \\ 37200-000 Lavras, Brazil  & 1050
Brussel, Belgium
 & 30.203 Cartagena, Spain \\ osnel@ufla.br & efjesper@vub.ac.be
  & Manuel.Ruiz@upct.es
\end{tabular}

\end{document}